\documentclass[sn-nature]{sn-jnl}

\usepackage{graphicx}%
\usepackage{multirow}%
\usepackage{amsmath,amssymb,amsfonts}%
\usepackage{amsthm}%
\usepackage{mathrsfs}%
\usepackage[title]{appendix}%
\usepackage{xcolor}%
\usepackage{textcomp}%
\usepackage{manyfoot}%
\usepackage{booktabs}%
\usepackage{algpseudocode}%
\usepackage{listings}%
\usepackage{bm}
\usepackage{enumitem}
\usepackage{accents}
\usepackage{booktabs}
\usepackage{dsfont}
\usepackage{longtable}
\usepackage[ruled]{algorithm2e}
\usepackage{hyperref}
\usepackage{subcaption}
\usepackage{float}
\usepackage{array}
\newcommand{\ubar}[1]{\underaccent{\bar}{#1}}

\newcommand{\ZZ}{\mathbb{Z}}
\newcommand{\RR}{\mathbb{R}}

\theoremstyle{thmstyleone}%
\theoremstyle{thmstyletwo}%

\theoremstyle{thmstylethree}%

\begin{document}

\title[Article Title]{Global Optimization: A Machine Learning Approach}


\author[1]{\fnm{Dimitris} \sur{Bertsimas}}\email{dbertsim@mit.edu}

\author[2]{\fnm{Georgios} \sur{Margaritis}}\email{geomar@mit.edu}

\affil[1]{\orgdiv{Sloan School of Management}, \orgname{Massachusetts Institute of Technology}, \orgaddress{\city{Cambridge}, \postcode{02139}, \state{MA}}}

\affil[2]{\orgdiv{Operations Research Center}, \orgname{Massachusetts Institute of Technology}, \orgaddress{\city{Cambridge}, \postcode{02139}, \state{MA}}}



\abstract{Many approaches for addressing Global Optimization problems typically rely on relaxations of nonlinear constraints over specific mathematical primitives.  This is restricting in applications with constraints that are black-box, implicit or consist of more general primitives. Trying to address such limitations, Bertsimas and Ozturk (2023) proposed OCTHaGOn as a way of solving black-box global optimization problems by approximating the nonlinear constraints using hyperplane-based Decision-Trees and then using those trees to construct a unified mixed integer optimization (MIO) approximation of the original problem. We provide extensions to this approach, by (i) approximating the original problem using other MIO-representable ML models besides Decision Trees, such as Gradient Boosted Trees, Multi Layer Perceptrons and Suport Vector Machines
	(ii) proposing adaptive sampling procedures for more accurate machine learning-based constraint approximations, 
	(iii) utilizing robust optimization to account for the uncertainty of the sample-dependent training of the ML models, 
	(iv) leveraging a family of relaxations to address the infeasibilities of the final MIO approximation.
	We then test the enhanced framework in 81 Global Optimization instances.
	We show improvements in solution feasibility and optimality in the majority of instances.
	We also compare against BARON, showing improved optimality gaps or solution times in 11 instances.}




\maketitle

\section{Introduction}
Global optimizers aim to solve problems of the following form:
\begin{equation}\label{eq:gop_intr}
	\begin{aligned}
		\min \quad & f(\bm{x})\\
		\textrm{s.t.} \quad &g_i(\bm{x})\leq 0,\: i\in \bar{I},\\
		&h_j(\bm{x})=0,\: j\in \bar{J},\\
		&\bm{x}\in\ZZ^m\times \RR^{n-m},
	\end{aligned}
\end{equation}
where $f,g_i,h_i$ represent the objective function,
the inequality constraints and the equality constraints
respectively.
The objective function and constraints may lack desirable
mathematical properties like linearity or convexity, and 
the decision variables may be continuous or integer.
\par
Most approaches in the Global Optimization literature,
attempt to solve Problem (\ref{eq:gop_intr}) by approximating it
with more tractable optimization forms.
For this purpose, they often use a combination of
gradient-based methods, outer approximations, relaxations
and Mixed Integer Optimization (MIO). 
For instance, the popular nonlinear Optimizer CONOPT
uses a gradient-based approach in its solution
process.
As noted in \cite{drud_conoptlarge-scale_1994},
CONOPT finds an initial feasible solution using heuristics,
performs a series of gradient descent iterations
and then confirms optimality via bound projections.
During the gradient-based part of the algorithm,
CONOPT linearizes the constraints and performs a series of
linear-search gradient-based iterations, while
preserving feasibility at each step.\par 
A different approach is the one detailed by \cite{horst_outer_1989},
which uses outer approximations.
This approach simplifies the problem by approximating
the constraints via linear and nonlinear cuts, while
preserving the initial feasible set of the problem.
Such approach can only be used with constraints that obey
a particular mathematical structure, such as
linearity of integer variables and convexity of
the nonlinear functions \cite{duran_outer-approximation_1986},
or concavity and bilinearity \cite{bergamini_improved_2008},
where the functions involved are amenable to efficient outer approximations.
Although such approaches are effective in some scenarios,
they have not found extensive use as they they are restricted
to certain classes of problems.\par 
Another approach is the one used
by the well-established commercial optimizer BARON, which
combines MIO with convex relaxations and outer approximations.
As detailed in \cite{ryoo_branch-and-reduce_1996-1}, BARON uses a
branch-and-reduce method, which partitions the domain
of the constraints and objective into subdomains,
and attempts to bound the decision variables in each
subdomain depending on the mathematical form of the constraints.
This approach produces a branch-and-bound tree that
offers guarantees of global optimality, which is analogous
to the branch-and-bound process used for solving MIO problems.
A similar approach is used by the popular solver ANTIGONE.
As described in \cite{misener_antigone_2014}, ANTIGONE
first reformulates the problem, detects specific mathematical
structure in the constraints and then uses convex underestimators
in conjunction with a branch-and-cut method
in order to recursively find the optimal solution.\par
Although all such approaches are very effective 
in dealing with certain types of global optimization problems,
they each have their own weaknesses.
For instance, gradient-based approaches depend on good
initial feasible solutions and cannot easily handle integer variables.
Also, approaches that leverage relaxations and outer
approximations are restricted to specific types
of nonlinearities, thus being ineffective against
more general problems.
Finally, approaches that use branching or MIO-based methods
inevitably suffer from the curse of dimensionality due
to their combinatorial nature.\par 
At the same time, for reasons of efficiency, 
many global optimizers only allow constraints that
use a subset of all the mathematical primitives.
For instance, the optimizer BARON can
only handle functions that involve exp($x$), ln($x$), 
$x^a$ and $b^x$ where $a$ and $b$ reals.
However, real-world problems may contain constraints with a much
richer set of primitives, such as trigonometric and
piecewise-discontinuous functions, where it is not
always possible to reduce the problem to a form
compatible with the global optimizer.
At the same time, such optimizers are generally
ineffective in cases where the constraints are
black-box and lack an analytical 
representation.
Such cases may, for instance, arise in 
simulation-based constraints.\par 
Motivated by such scenarios, 
\cite{bertsimas_global_2022} proposed OCTHaGOn, 
a Global Optimization framework that attempts 
to solve general and black-box Global Optimization 
problems by combining Machine Learning with Mixed 
Integer Optimization (MIO). 
The framework can, in theory, be used in problems that involve any type
of constraints, such as convex, non-convex, implicit and constraints that
involve very general primitives.
The only requirement of the method is that the decision variables
involved in nonlinear inequalities are bounded.
Under this assumption, OCTHaGOn first samples the nonlinear constraints 
for feasibility and trains a hyperplane-based 
decision tree (OCT-H, \cite{bertsimas_optimal_2017},\cite{bertsimas_machine_2019}) on those samples.
It does that for every non-linear constraint 
and then uses the MIO-representation of the 
resulting decision trees to construct an MIO
approximation of the original problem. 
In then solves approximating MIO and uses 
a local-search Projected Gradient Descent 
method to improve the solution in 
terms of both feasibility and optimality.\par 
Despite the generality of this approach, the framework
sometimes yields infeasible or suboptimal solutions
due to being approximate in nature.
Hence, in this paper, we provide extensions
to OCTHaGOn in order to improve both feasibility and optimality 
of the generated solutions. 
We enhance the framework by leveraging
adaptive sampling, robust optimization, constraint relaxations and a variety 
of MI-representable ML models to create better approximations of the original problem.
We test the enhanced framework in a number of Global Optimization benchmarks,
and we show improved optimality gaps and solution times in the majority of test instances.
We also compare against the commercial optimizer BARON, showing improvements in a number of
instances.

\subsection{Machine Learning in Optimization}
The interplay between Machine Learning and  Optimization has
been an active area recently. 
\cite{sun_survey_2019} provide a survey of different optimization
methods from an ML standpoint.
In another survey, \cite{gambella_optimization_2021} demonstrate
how fundamental ML tasks, such as classification, regression and 
clustering can be formulated as optimization problems.
Additionally, 
a lot of emphasis has been placed on using ML for aiding 
the formulation and solution of optimization problems.
For instance, ML has been used to assist with the solution
of difficult and large-scale optimization formulations
(e.g. \cite{fischetti_machine_2019}, \cite{abbasi_predicting_2020}),
and it has also been used to speed-up tree search search \cite{hottung_deep_2020}
and guide branching in MIO problems \cite{khalil_learning_2016}.\par 
Out of the various works on the intersection of ML and Optimization,
a very relevant one to our approach is the framework of Constraint
Learning \cite{maragno_mixed-integer_2021}.
In their work, the authors propose a systematic way of
(1) learning functions
from data through MIO representable ML models, (2)
embedding those learned functions as constraints in an MIO formulation
and (3) 
solving the resulting MIO in a way that these data-driven constraints
are satisfied.
The authors then use this framework to design food baskets
for the World Food Program, with ingredients compositions that satisfy
data-driven palatability constraints.
In a later survey, \cite{fajemisin_optimization_2022} demonstrate
various applications of this framework,
and then showcase preexisting works that utilize
its different parts.\par 
The framework OCTHaGOn \cite{bertsimas_global_2022} is an adaptation
of Constraint Learning for solving 
Global Optimization problems.
OCTHaGOn samples the non-linear constraints and learns the non-linearities
using MIO representable decision trees.
The type of decision trees used in the approximation are called Optimal Classification/Regression
Trees (OCT-Hs/ORTs, \cite{bertsimas_optimal_2017},\cite{bertsimas_machine_2019})
The trained decision trees are 
embedded as constraints into a Mixed Integer Optimization  (MIO)   formulation that approximates the original problem.
After a solution to the MIO is obtained,
a local improvement is performed to address for infeasibility
and suboptimality which stem from the approximation errors.\par 

\subsection{Contributions}
Our approach for solving Global Optimization Problems is an improvement of
OCTHaGOn with the following key differences and enhancements: 
\begin{enumerate}[label={(\roman*)}]
	\item \textbf{Different ML models:}
	OCTHaGOn only uses hyperplane-based decision trees (OCT-H) as constraint approximators, 
	but we extend the approach to other types of MIO-representable
	ML models, such as Support Vector Machines, Gradient Boosted Trees, and Multi-layer
	Perceptrons.
	We also employ a method for selecting which of those models to use for
	which constraint.
	\item \textbf{Adaptive Sampling:}
	We propose an adaptive sampling procedure that helps generate
	high-quality samples of the nonlinear constraints for more accurate constraint approximations.
	This procedure is adaptive and uses hyperplane-based decision trees
	to over-sample hard-to-approximate areas of the nonlinear constraints.
	
	\item \textbf{MI Relaxations:}
	The ML-based MIO approximation of the original problem can sometimes be infeasible.
	To eliminate such infeasibilities, we relax constraints
	of the resulting MIO in a way that makes the MIO approximation
	feasible if the original problem is also feasible.
	\item \textbf{Robust Optimization}:
	We utilize Robust Optimization to account
	for the uncertainty of the sample-dependent training of the ML models. In particular, we
	attempt to model the uncertainty in the trained parameters of the ML approximators,
	and we then use Robust Optimization to correct for this uncertainty. 
	
\end{enumerate}
We test the framework on $81$ Global Optimization instances, with $77$ of those
being part the standard benchmarking library MINLPLib.
Our results show that in the majority of instances, our enhancements 
improve the optimality gaps or solutions times of OCTHaGOn.
We also show that in $11$ out of the $81$ instances, the enhanced framework 
provides better or faster solutions than BARON.\par  
Overall, we show that despite its approximate nature, 
the enhanced framework, is a promising method
in finding globally optimal solutions.
The framework can be applied in problems with convex, non-convex
and even implicit and black-box constraints, with the only assumption
that the user needs to specify bounds for the decision variables.
Hence, due to its generality the method can be used in problems that
are incompatible with traditional global optimizers such as BARON and ANTIGONE.

%
%
%

\section{OCTHaGOn}
\label{section:octh}
Since our method is an extension of OCTHaGOn, we will use this section to
present the steps involved in the OCTHaGOn framework.
OCHaGOn \cite{bertsimas_global_2022} attempts to solve
Global Optimization problems
using a combination of ML and MIO.
It relies on creating decision-tree based approximations of the nonlinear constraints,
embedding those approximations in an approximating MIO, and then solving and improving the solution
resulting from that MIO to obtain a high-quality solution of the original problem.\par 
Due to its nature, OCTHaGOn can also be used to solve problems with
very general black-box, implicit or even data-driven constraints.
The only requirement is that the user of the framework supplies bounds 
for the variables involved in the nonlinear constraints.
If such bounds are not provided, the quality of the solution returned
by the framework can be adversely affected. 
\par 
\subsection{Algorithmic process}
OCTHaGOn addresses problems of the following form:
\begin{equation}\label{eq:gop}
	\begin{aligned}
		\min \quad & f(\bm{x})\\
		\textrm{s.t.} \quad &g_i(\bm{x})\leq 0,\: i\in \bar{I},\\
		&h_j(\bm{x})=0,\: j\in \bar{J},\\
		&\bm{x}\in\ZZ^m\times \RR^{n-m},
	\end{aligned}
\end{equation}
where the functions $f,g_i,h_j$ can be non-convex.



OCTHaGOn attempts to create and solve an MIO approximation of Problem (\ref{eq:gop}) in the following steps:
\begin{enumerate}
	\item \textbf{Standard form problem generation:}
	In order to create the MIO approximation of Problem (\ref{eq:gop}), 
	we first separate the linear from
	the non-linear constraints.
	We also find the constraints that define explicit bounds $[\underbar{x}_k,\overline{x}_k]$
	for every optimization variable $x_k$ that is involved in a nonlinear constraint.
	If a variable $x_k$ does not  have an explicit bound (either from above or below),
	we attempt to compute that bound by minimizing/maximizing that variable 
	subject to the linear constraints of (\ref{eq:gop}).
	Note that having explicit variable bounds is essential for 
	 Step 2 of the methodology.
	Then, after determining the variable bounds, we rewrite our original problem into
	the following form:
	\begin{equation}\label{eq:gop_std}
		\begin{aligned}
			\min_{\bm{x}} \quad & f(\bm{x})\\
			\textrm{s.t.} \quad &g_i(\bm{x})\leq 0,\: i\in {I},\\
			&h_j(\bm{x})=0,\: j\in {J},\\
			&\bm{Ax}\geq \bm{b},\: \bm{Cx}=\bm{d},\\
			&x_k\in[\underbar{x}_k,\overline{x}_k],\: k\in[n].
		\end{aligned}
	\end{equation}
	The linear constraints of this resulting formulation are directly passed to the MIO 
	approximation, while the non-linear ones are approximated using Steps 2-6.
	\item \textbf{Sample and evaluate nonlinear constraints:}
	For each nonlinear constraint $g_i(\bm{x})\leq 0$, we generate samples
	of the form $D_i=\{(\bm{\widetilde x_k},\tilde{y}_k)\}^n_{k=1}$
	where $\tilde{y}_k=\mathds{1}\{g(\bm{\widetilde x_k})\leq 0\}$ and 
	$\mathds{1}\{\cdot \}$ is the indicator function.
	In order to chooses the points $\bm{\widetilde x_1},\dots,\bm{\widetilde x_k}$
	where the samples are evaluated,  the framework uses the following
	sampling methods: Boundary Sampling, Latin Hypercube sampling (LH sampling) 
	and kNN Quasi-Newton sampling.
	Boundary sampling samples the corners of the hyper-rectangle
	that is formed by the bounds of the decision variables (i.e., for each
	decision variable $x_j$ of the constraint, we have that $x_j\in[\ubar{x}_j,\bar{x}_j]$
	after computing bounds in Step 1).
	LH sampling, due to its space-filling characteristics, is then used to obtain constraint
	samples from the interior
	of the decision-variable hype-rectangle (i.e., each decision variable $x_j$ is sampled
	in the range $[\ubar{x}_j,\bar{x}_j]$).
	Finally, kNN Quasi-Newton sampling is proposed by the authors
	as a way to sample the boundaries of the constraints, given
	the already generated samples from the previous steps.
	
	\item \textbf{Train decision tree:}
	Hyperplane based decision trees (OCT-H, \cite{bertsimas_optimal_2017},\cite{bertsimas_machine_2019}) are trained on the datasets
	$D_i$ to approximate the feasibility space
	of the nonlinear constraint $g_i(\bm{x})\leq 0$.
	This yields decision tree aproximators $\bar{c}_i(\bm{x})$ such that 
	$\bar{c}_i(\bm{x})\simeq\mathds{1}\{g_i(\bm{x})\leq 0\}$.
	If we have a non-linear objective, then this objective is approximated
	with regression trees (ORT-H, \cite{bertsimas_optimal_2017},\cite{bertsimas_machine_2019})
	instead of classification trees.
	
	\item \textbf{Generate MIO approximation:}
	The trained decision trees $\bar{c}_i(\bm{x})$ are transformed into an efficient
	MIO representation by extracting the corresponding hyperplane splits and
	then using disjunction formulations.
	The resulting representations are then embedded as constraints of the form $\bar{c}_i(\bm{x})=1$
	into the approximation MIO of Problem (\ref{eq:gop}).
	
	\item \textbf{Solve MIO approximation:}
	The MIO approximation of Problem (\ref{eq:gop}) is solved using 
	commercial solvers, such as IBM CPLEX and Gurobi.
	
	\item \textbf{Check and improve solution:}
	The MIO is just an approximation of the original global optimization
	problem, so the solution of Step 5 can be near-optimal and near-feasible.
	In order to account for this, the framework first measures the feasibility
	of the solution with respect to the nonlinear constraints, then
	it computes the gradient of the nonlinear objective and constraints
	(i.e., using automatic differentiation), and finally it performs a sequence of 
	projected gradient descent steps to improve feasibility and optimality
	of the solution.
\end{enumerate}
The authors of OCTHaGOn test the framework against
various benchmarks of the MINLP library, and demonstrate
that in many cases, it works well compared to traditional
global optimizers such as CONOPT, IPOPT and BARON.

\section{Enhancements}
In this section, we propose a variety of improvements
on top of OCTHaGOn which improve the quality
of the solutions generated by the framework.
We first provide a brief outline of the solution steps
of the enhanced framework.
For each step, we also indicate whether the step is the same
as OCTHaGOn, whether it's an enhanced version of an OCTHaGOn
step, or whether it's a new step.
The outline of the different steps are shown below:

\begin{enumerate}
	\item \textbf{Standard form problem generation:}
	We separate the linear from the nonlinear constraints and
	attempt to compute bounds for the variables involved in the nonlinear constraints.
	Then, we initialize the MIO approximation of the original problem with the linear
	and bound constraints.
	The non-linear constraints will be approximated and included in this MIO approximation
	in Steps 2-4. 
	\item \textbf{Sample and evaluate nonlinear constraints (Enhanced Step):}
	We sample the nonlinear constraints  in the following   steps:
	(i) Boundary sampling, (ii) Latin Hypercube Sampling, (iii) KNN Quasi-Newton sampling
	and (iv) OCT-based adaptive sampling.
	The first $3$ sampling steps are the same as the ones used in OCTHaGOn and 
	are used to obtain an initial set of samples.
	Then, we use OCT-based adaptive sampling, which is part of our enhancements,
	to adaptively obtain high-quality samples in areas where the nonlinear 
	constraints are not approximated well by the ML learners (e.g. near the constraint
	boundaries).
	
	\item \textbf{ML model training (Enhanced Step):}
	For each nonlinear constraint, we train a series of MIO-representable ML models
	(namely OCT-Hs, MLPs, SVMs and GBMs) on the samples generated
	during Step $2$.
	We measure the accuracy of those models on a test set, and
	for each nonlinear constraint, we pick the model that demonstrates
	the best predictive performance.
	We then use the model for approximating the corresponding nonlinear
	constraint.
	We follow a similar process to approximate the non-linear objective
	(if such objective is present), but in this case, we use regression instead of
	classification models.
	The enhancement in this step is the training of multiple ML models (MLPs, SVMs, GBMs),
	besides the decision trees (OCT-Hs) which are the only models used by OCTHaGOn.
	\item \textbf{Generate MIO approximation (Enhanced Step):}
	For each one of the nonlinear constraints, we take the corresponding ML approximator
	from the last step and we represent it using an MIO formulation.
	Then, we embed all the MIO approximations into a unified MIO formulation.
	In this formulation, we also include the linear constraints (and objectives)
	of the original problem.
	The enhancement in this step is that besides the decision trees (OCT-Hs) used
	by OCTHaGOn, we also formulate and represent the additional
	ML models (MLPs, SVMs and GBMs) using a MIO formulation.
	\item \textbf{Introduce Robustness (New Step):}
	In order to account for the uncertainty in the trained
	parameters of the ML constraint approximators, we utilize 
	Robust Optimization.
	We assume that the trained parameters of SVMs, GBMs and OCT-Hs lie in $p$-norm
	uncertainty sets, and we modify the MIO representation of the learners
	by considering the corresponding robust counterparts.
	The level of robustness is controlled by the parameter $\rho$,
	where a value of $\rho=0$ means that we do not account for robustness,
	while higher values of $\rho$ lead to more robust and more conservative
	solutions.
	
	\item \textbf{Constraint Relaxation (New Step):}
	In certain cases, the MIO approximation of the original
	problem is infeasible, due to the inexact nature of the constraint
	approximators. 
	To account for this issue, we relax the MI constraints,
	while introducing a relaxation penalty into the objective. 
	
	\item \textbf{Solve MIO approximation:}
	We solve the MIO approximation of the original problem
	using commercial solvers, such as IBM CPLEX and Gurobi.
	
	\item \textbf{Improve solution (Enhanced Step):}
	Since the MIO solution can be near-optimal and near-feasible, 
	we perform a series of Projected Gradient Descent iterations in
	order to improve feasibility and optimality of the solution.
	In this PGD step, we use the same procedure as in OCTHaGOn,  but
	we also conditionally use Momentum
	in order to reduce the chances of landing in local optima.
	
%
	
\end{enumerate}
\noindent
We next elaborate on each one of the proposed enhancements.

\subsection{Enhancement 1: ML Model Training}
\label{section:ml-models}
One of the key aspects of OCTHaGOn framework
is the approximation of constraints using MIO-representable ML
models and the incorporation of those models into
a common MIO.
For this approximation, OCTHaGOn uses hyperplane-based
decision trees (OCT-H and ORT-H, \cite{bertsimas_optimal_2017},\cite{bertsimas_machine_2019}) which are
improved generalizations of Classification and Regression Trees (CART).\par 
However, OCT-Hs are clearly not the only models that can
be used to for constraint approximations.
In fact, there are many other types of MIO-representable
ML models, such as Linear Regression models, Support Vector Machines (SVMs),
Random Forests (RFs), Gradient Boosting Machines (GBMs) and
multi-layer perceptrons with ReLU activations (MLPs).\par 
Hence, given that different types of models may have different 
strengths and weaknesses, we chose to incorporate a wider
family of models (i.e., SVMs, GBMs and MLPs), 
while also employing a cross-validation procedure
to select which model type to use for each constraint.
Our approach for incorporating different types of models
and selecting model types was based on 
\cite{maragno_mixed-integer_2021} and will be analyzed below.\par 
Before proceeding into the analysis of each model type, we  
note that all the ML models we use,
support both classification and regression.
In particular, for nonlinear objectives $f(\bm{x})$,
we create samples $D_R=\{(\bm{\widetilde x_k},f(\bm{\widetilde x_k}))\}^N_{k=1}$
and we train a regressor on those samples,
whereas for nonlinear constraints $g(\bm{x})\leq 0$ we create
samples  $D_C=\{(\bm{\widetilde x_k},\mathds{1}\{g(\bm{\widetilde x_k})\leq 0\})\}^N_{k=1}$
and we train a classifier on those samples.
Hence, when representing the different ML models
using MIO constraints, we will describe separately the case 
for classifiers and regressors.


\par 
\begin{itemize}
	\item \textbf{Support Vector Machines}:
	Support vector machines are ML models that use a suitable
	hyperplane to make predictions, either for classification
	\cite{cortes_support-vector_1995} or regression \cite{drucker_support_1996}.
	For the purpose of this work, we are only 
	considering the case of linear SVMs and we are not
	using any type of kernel, since our goal is to 
	generate linear and tractable MIO approximations of the constraints.\par 
	In the case of the linear support vector regression (SVR),
	we follow \cite{maragno_mixed-integer_2021,drucker_support_1996} to fit
	a linear function to the regression data $D_R$.
	The procedure is very similar to linear regression, with the difference
	that the loss function only penalizes residuals greater than a threshold $\epsilon$.
	As with linear regression, the resulting model is a linear function of the input
	with learned coefficents $\bm{\beta},\beta_0$ and can be embedded to the optimization
	as follows:
	\begin{equation}
		y_{SVR}=\beta_0+\bm{\beta}^T\bm{x}.
	\end{equation}
	In a classification setting,
	support vector classifiers (SVC) are trained
	to find a hyperplane that best separates positive and negative samples \cite{cortes_support-vector_1995}.
	After the SVC is trained on the classification data $D_C$, 
	we can use the model parameters $\beta_0$ and $\bm{\beta}$,
	to make binary predictions as follows: 
	\begin{equation}
		y_{SVC} = \mathds{1}(\beta_0+\bm{\beta}^T\bm{x}\geq 0).
	\end{equation} 
	Note that when training the SVC, the binary labels are first converted to $\{-1,1\}$
	from the binary form $\{0,1\}$ used in dataset $D_C$.
	Then, we can approximate the nonlinear constraint $g(\bm{x})\leq 0$
	with the linear constraint $\beta_0+\bm{\beta}^T\bm{x}\geq 0$.
	
	\item \textbf{Decision Trees}:
	Decision trees are ML models that partition the observation
	space into disjoint leaves after a sequence of splits.
	Due to their inherent interpretability and their ability to model
	non-linear relationships in data, they have been used a lot in practice.
	\cite{breiman_classification_1984} first introduced Classification and Regression
	Trees (CART), a framework that greedily partitions the observation
	space using parallel splits (i.e., splitting by one feature at a time).
	Ever since, decision tree models have been revised and extended.
	\cite{bertsimas_optimal_2017} proposed a generalization
	of decision trees, where the decision-tree generation
	problem is formulated as an MIO and solved to near optimality.
	This approach also supports splitting over multiple features
	at a time by using hyperplane based splits (instead of parallel splits).\par 
	For the purpose of this work, we will model the generalized version
	of decision trees which includes hyperplane splits, since this version can also
	capture the standard CART trees with parallel splits.
	In particular, a decision tree can be modeled as a binary tree
	with non-terminal (intermediate) and terminal (leaf) nodes.
	Each non-terminal node $N_i$ represents a split of the form $\bm{a}^T_i\bm{x}\leq b_i$,
	and each leaf node $L_i$ is associated with a prediction $p_i$,
	where $\bm{a}_i,b_i,p_i$ are
	model parameters learned through training.
	Then, given an input vector $\bm{x}$, we can make predictions as follows:
	starting from the root $N_1$, we check whether $\bm{a}^T_1\bm{x}\leq b_1$.
	If the condition is satisfied, then we proceed to the left child
	of the root.
	Otherwise, we proceed to the right child.
	We recursively repeat this process by checking the split condition
	of the nodes and continuing traversing the tree until we reach a leaf node.
	When a leaf node $L_i$ is reached, we output the prediction $p_i$.
	If we now use $L(L_i)$ and $R(L_i)$ to denote the set of non-terminal nodes
	for which leaf $L_i$ is contained in their left and right subtree respectively,
	then each leaf $L_i$ is described by the following polyhedron:
	\begin{equation}
		\mathcal{P}_i = \Big(\bigcap_{j\in L(L_i)}\{\bm{x}\in\RR^n:\bm{a}_j^T\bm{x}\leq b_j\}\Big)\bigcap\Big(\bigcap_{j\in R(L_i)}\{\bm{x}\in\RR^n:\bm{a}_j^T\bm{x}> b_j\}\Big).
	\end{equation}
	Then, the decision tree predicts $p_i$ if and only if $\bm{x}\in\mathcal{P}_i$.
	Note here that the polyhedra $\mathcal{P}_1,\mathcal{P}_2,\dots$ are disjoint and their union  is MIO representable.
	Additionally,  this representation can also be used
	to describe decision trees with parallel splits (i.e., CART).
	In particular, if all the vectors $\bm{a}_i$ are binary with $\sum_j{a}_i^{(j)}=1$,
	then all the splits in the tree are done in a single feature at a time,
	and the polyhedrons $\mathcal{P}_i$ are hyper-rectangles.
	\par 
	Let's now use $\mathcal{L}$ to denote the set of leaves.
	If we want to embed the tree predictions into our MIO under a regression setting,
	then we can use the following set of constraints to do so:
	\begin{equation}\label{eq:dt}
		\begin{aligned}
			\sum_{i=1}^{|\mathcal{L}|}z_ip_i&=y_{DT},\\
			\sum_{i=1}^{|\mathcal{L}|}z_i&=1,\\
			\bm{a}_j^T\bm{x}&\leq b_j + M(1-z_i), \quad \forall i\in\mathcal{L},j\in L(L_i),\\
			\bm{a}_j^T\bm{x}&\geq  b_j - M(1-z_i)+\epsilon, \quad \forall i\in\mathcal{L},j\in R(L_i),
		\end{aligned}
	\end{equation}
	where $y_{DT}$ represents the output of the decision tree and
	$\bm{z}\in\{0,1\}^{|\mathcal{L}|}$ is a vector of binary variables.\par 
	On the other hand, if we have a classification decision tree trained
	on the dataset $D_C=\{(\bm{\widetilde x_k},\mathds{1}\{g(\bm{\widetilde x_k})\leq 0\})\}^N_{k=1}$,
	then we can require that we are in a feasible leaf of that tree
	by imposing the following additional constraint:
	\begin{equation}
		y_{DT}\geq 0.5.
	\end{equation}
	
	\item \textbf{Gradient Boosting Machines}:
	Gradient Boosted Machines (GBM) can be represented as 
	ensembles of base-learners, where each learner is trained
	sequentially to correct the mistakes of the previous learner.
	After training, if we use $y_i$ to denote the prediction of the 
	$i$-th learner, then the output of the GBM model can be described
	as:
	
	\begin{equation}
		y=\sum_{i=1}^n a_iy_i,
	\end{equation}
	where $a_i$ are weights associated with each learner.
	In theory, many different types of models can be used
	as based-learners for GBM, but we will limit our focus
	to Decision Trees, which are used most often in practice.
	In that case, we can embed the predictions of GBMs
	into the MIO by adding the following constraint:
	\begin{equation}\label{eq:gbm}
		y_{GBM}=\sum_{i=1}^n a_i\:y^{(i)}_{DT},
	\end{equation}
	where $y^{(i)}_{DT}$ represents the predictions
	of the $i$-th tree base-learner.
	Hence, the prediction of the GBM is a weighted
	average of the predictions of the individual
	decision trees $y^{(i)}_{DT}$, where
	the MIO representation of $y^{(i)}_{DT}$ is the 
	one we described in Eq. (\ref{eq:dt}). 
	Then, in order to use to approximate the constraint $g(\bm{x})\leq 0$
	using a GBM, we first train the GBM on the classification
	dataset $D_C=\{(\bm{\widetilde x_k},\mathds{1}\{g(\bm{\widetilde x_k})\leq 0\})\}^N_{k=1}$,
	we represent the output of the GBM using the constraint
	(\ref{eq:gbm}) and then we require that $y_{GBM} \geq 0.5$.
	
	\item \textbf{Neural Networks}: A powerful and popular class of Neural Network models consists of
	multi-layer perceptrons (MLPs) with rectified linear unit activations (ReLU).
	Such networks consist of an input layer, $L-2$ hidden layers with
	ReLU activations and an output layer.
	The repeated use of the nonlinear ReLU operator ($\max\{0,x\}$) in nodes of such networks
	allows the network to approximate a variety of nonlinear functions very well, leading
	to a class of models with very high representational power.\par 
	
	Additionally, ReLU-based MLPs are also MIO-representable
	through a linear big-M formulation \cite{grimstad_relu_2019}.
	In particular, following \cite{maragno_mixed-integer_2021}, if we define $N^l$ 
	as the set of nodes of $l$-th hidden layer of the network,
	then the value $v_i^l$ of each node $i\in N^l$ is
	calculated as a weighted average of the node
	values of the previous layer.
	The result is passed through a ReLU activation yielding the following formula:
	\begin{equation}
		v_i^l=\max\Biggl\{0,\quad \beta_{i0}^l+\sum_{j\in N^{l-1}}\beta_{ij}^lv_j^{l-1}\Biggr\},
	\end{equation} 
	\noindent
	where $\beta_{i}^l$ is the vector of coefficients for
	node $i$ in layer $l$ retrieved after training.
	The benefit of such models is that due to
	their nonlinearities, they can be used to model very complex
	nonlinear constraints and objectives compared to other types of models.
	\par 
	For our use-case, we will only consider
	MLPs with $1$ output neuron.
	In particular, for the regression task, we train an MLP regressor on
	the dataset $D_R=\{(\bm{\widetilde x_k},f(\bm{\widetilde x_k}))\}^n_{k=1}$
	using Adam Optimizer with Mean Squared Error (MSE) loss. 
	Then, we can use the following constraints to represent the output
	of the MLP regressor:
	\begin{equation}\label{eq:mlp}
		\begin{aligned}
			y_{MLP}=&\beta_{00}^L+\sum_{j\in N^{L-1}}\beta_{0j}^Lv_j^{L-1},\\
			u_i^l&\geq \beta_{i0}^l+\sum_{j\in N^{l-1}}\beta_{ij}^lv_j^{l-1},\quad &\forall l=\{2,\dots,L-1\},i\in N^l,\\
			u_i^l&\leq \beta_{i0}^l+\sum_{j\in N^{l-1}}\beta_{ij}^lv_j^{l-1}-M(1-z_{il}),\quad& \forall l=\{2,\dots,L-1\},i\in N^l,\\
			u_i^l&\leq Mz_{il},\quad &\forall l=\{2,\dots,L-1\},i\in N^l,\\
			u_i^l&\geq 0,\quad &\forall l=\{2,\dots,L-1\},i\in N^l,\\
			u_i^1 &= x_i,\quad &\forall i\in[n],\\
			z_{il}&\in\{0,1\},
		\end{aligned}
	\end{equation}
	where $x_i$ is the input variable and $y_{MLP}$ is the model's output,
	which is an approximator of the objective $f(\bm{x})$.
	Here, we rely on a big-M formulation to model the ReLU activations,
	which can be tightened through an appropriate choice
	of $M$ \cite{maragno_mixed-integer_2021}.\par 
	For classification, we train an MLP on
	the dataset $D_C=\{(\bm{\widetilde x_k},\mathds{1}\{g(\bm{\widetilde x_k})\leq 0\})\}^N_{k=1}$ using binary cross-entropy loss
	and a sigmoid activation on the output node.
	Then, we approximate the constraint $g(\bm{x})\leq 0$
	using the formulation (\ref{eq:mlp}),
	but with the additional constraint $y_{MLP}\geq 0$.
	This new constraint is used to ensure that the output logit
	is positive, which corresponds to an output probability
	greater than $0.5$. 
	Note that we do not need to explicitly model the nonlinear sigmoid activation in
	our formlulation, we can simply disregard the output sigmoid node and require that
	its input (logit) is positive.\\
	
\end{itemize}
\leavevmode \noindent
In this work, we can use
all the ML models described above to approximate
a non-linear objective or constraint.
Choosing which model type to use for which 
nonlinear function is an important task, since this choice heavily
determines the quality of the resulting MIO approximation.
For this reason, we employ a cross-validation 
procedure to determine the best model type for each nonlinear
function.
Following \cite{maragno_mixed-integer_2021}, we train
all model types on samples of the nonlinear function,
we measure the accuracy (for classifiers) and $R^2$ (for
regressors) in a test set of unseen function samples,
and we choose the ML model that has the best performance.\par 
Then, we approximate every non-linear constraint $g_i(\bm{x})\leq 0$ of Problem (\ref{eq:gop})
with the following constraint:
\begin{equation}
	y^{(i)}_{ML}(\bm{x})\geq a_i,
\end{equation}
where $y^{(i)}_{ML}(\bm{x})$ is the MIO representation of the best-performing
model for the constraint $g_i(\bm{x})\leq 0$ and $a_i$ is a fixed threshold
that depends on the model type.
For each model type, the corresponding MIO approximations $y^{(i)}_{ML}(\bm{x})$
and thresholds $a_i$ are the ones we described in the previous paragraphs.
For instance, for MLPs we have that $y^{(i)}_{ML}(\bm{x})=y_{MLP}$ and
$a_i=0$, and for GBMs we have $y^{(i)}_{ML}(\bm{x})=y_{GBM}$ and $a_i=0.5$. 
Similarly, we approximate a non-linear objective $f(\bm{x})$ using the MIO representation
of the best performing regression model (i.e., $f(\bm{x})\simeq y_{ML\_REGR}(\bm{x})$).\par 
More concretely, if we have the following optimization
problem:
\begin{equation}
	\begin{aligned}
		\min_{\bm{x}\in\RR^n} \quad & f(\bm{x})\\
		\textrm{s.t.} \quad &g_i(\bm{x})\leq 0,\: i\in I,\\
		&\bm{Cx}\leq  \bm{d},
	\end{aligned}
\end{equation}
where $f(\bm{x})$ and $g_i(\bm{x})$ are nonlinear functions,
then the MIO approximation of that problem has the following form:
\begin{equation}\label{eq:mio}
	\begin{aligned}
		\min_{\bm{x}\in\RR^n} \quad & y_{ML\_REGR}(\bm{x})\\
		\textrm{s.t.} \quad &y^{(i)}_{ML}(\bm{x})\geq a_i,\: i\in I,\\
		&\bm{Cx}\leq  \bm{d},
	\end{aligned}
\end{equation}
where $y_{ML\_REGR}(\bm{x})$ and $y^{(i)}_{ML}(\bm{x})$ 
are MIO representable functions.\bigskip


\subsection{Enhancement 2: Sampling}
\label{section:gop_sampling}
In order to have good approximations of the
nonlinear constraints, we first need to have high-quality
samples of those constraints.
Our samples need to capture 
the non-linear and non-convex feasible regions of the constraints,
which means that a static sampling procedure
is not enough for our task.
In particular, the sample-generation process
should not be agnostic to the sampled constraint,
but it should instead adapt to the landscape of
the function we are attempting to sample.\par 
As mentioned in Section \ref{section:octh}, the sampling 
procedure of OCTHaGOn involves three steps:
Boundary Sampling, LH sampling and kNN Quasi-Newton
sampling.
Both Boundary Sampling and LH sampling are static
sampling procedures, in the sense that they sample all functions in
the same manner, without taking into account the
characteristics of those functions to determine where to sample.
On the other hand, kNN Quasi-Newton sampling
is a more adaptive sampling procedure, given
that it uses the function output to sample 
the constraints on the boundaries of the feasibility region.\par 
For our approach, we initially use the same three steps as
OCTHaGOn in order to generate a first set of samples.
Then, we propose an adaptive sampling method in order
to generate even more fine-grained samples for the constraints.
Our method is called OCT Sampling and attempts 
to iteratively resample the constraint
in areas where the learners cannot approximate well.
The method is described in detail below.\par 
\textit{OCT Sampling}:
Let's assume that we want to approximate the constraint
$g(\bm{x})\leq 0$ and we already 
have samples $D=\{(\bm{\widetilde x_k},\mathds{1}\{g(\bm{\widetilde x_k})\leq 0\})\}^N_{k=1}$
for the feasibility of that constraint.
The goal of OCT sampling is to resample parts
of the constraint which are difficult
to approximate.
In order to do that, we will use hyperplane-based
decision trees (OCT-H, \cite{bertsimas_optimal_2017},\cite{bertsimas_machine_2019})
which are one of the types of learners we use for constraint
approximations.\par 
The procedure we follow is the following:
\begin{enumerate}
	\item \textbf{Learner training}:
	We generate $K$ random subsets $D_1,\dots,D_K\subset D$
	of the dataset $D$, with a fixed size $|D_i|=C,\:\forall i\in[K]$.
	We then train one OCT-H learner $T_i$ on each dataset
	$D_i$.
	\item \textbf{Identify ambiguous samples}:
	In this step, we identify points $\bm{x}_1,\bm{x}_2,\dots$
	for which there is a high prediction discordance
	between the learners $T_1,\dots,T_K$.
	The goal of identifying such points is that 
	those points are indicators of areas where there
	is a poor generalization from our learners,
	and thus areas that are worth resampling.
	In order to find such points, let's use $T_i(\bm{x})$ to denote the binary
	prediction of the $i$-th learner given an input vector $\bm{x}$.
	We then define the following quantities:
	\begin{equation}
		\begin{aligned}
			P(\bm{x})&=|\{i\in[K]:T_i(\bm{x})=1\}|,\\
			N(\bm{x})&=|\{i\in[K]:T_i(\bm{x})=0\}|.\\
		\end{aligned}
	\end{equation}
	Then, we find a subset of points of $D$ for which there is a high
	discordance between the predictions of $T_1,\dots,T_K$.
	This set of points $S$ is defined as follows.:
	\begin{equation}
		S=\{\bm{x}\in D:|P(\bm{x})-N(\bm{x})|\leq K\tau \},
	\end{equation}
	where $\tau\in[0,1]$ is a threshold that determines
	the target level of predictor discordance.
	The higher the value of $\tau$, 
	the less the level of discordance is needed
	for a point $\bm{x}$ to be included in $S$.
	
	
	\item \textbf{Identify ambiguous polyhedra}:
	The goal of this step is to identify polyhedral
	regions where the points of $S$ reside.
	As we discussed in Section \ref{section:ml-models},
	in a hyperplane-based decision tree $T_i$,
	every leaf $L^{(i)}_j\in\mathcal{L}_i$ is represented by
	a polyhedron $\mathcal{P}_j^{(i)}$, while the polyhedra
	of the different leaves are disjoint.
	Hence, the decision tree $T_i$ assigns
	every point $\bm{x}\in \RR^n$
	to exactly $1$ leaf-polyhedron, which we will denote
	as $\mathcal{P}^{(i)}(\bm{x})$.
	Then, for every point $\bm{x}\in S$,
	we take the intersection $\mathcal{P}(\bm{x})$ of those leaf polyhedras
	of the different trees $T_1,\dots,T_K$, where:
	\begin{equation}
		\mathcal{P}(\bm{x})=\bigcap_{i=1}^K\mathcal{P}^{(i)}(\bm{x}).
	\end{equation}
	Note that $\mathcal{P}(\bm{x})$ is a polyhedron as
	an intersection of polyhedras.
	We then create the set of polyhedras 
	$\mathcal{C}=\{\mathcal{P}(\bm{x}):\bm{x}\in S\}$,
	which has the following property:
	\begin{equation}
		|P(\bm{x})-N(\bm{x})|\leq K\tau, \quad \forall \bm{x}\in\mathcal{P},\: \mathcal{P}\in\mathcal{C}.
	\end{equation}
	This means that if we pick any point $\bm{x}$ from any
	polyhedron $\mathcal{P}\in\mathcal{C}$, 
	then the predictions of the trees $T_1,\dots T_K$ for that point are
	guaranteed to have a particular level of disagreement.
	
	\item \textbf{Sample ambiguous polyhedra}:
	The goal of this step is to generate samples 
	from the interior of the polyhedra $\mathcal{P}\in\mathcal{C}$.
	As we mentioned before, those polyhedra represent
	areas where there is a high level of disagreement 
	between the predictors $T_1,\dots,T_K$.
	In order to sample those polyhedra, we will
	use the hit-and-run algorithm \cite{smith_efficient_1984}
	which is a Markov Chain Monte Carlo (MCMC) method for
	generating samples in the interior of a convex body.
	In particular, given a non-empty polyhedron $\mathcal{P}=\{\bm{x}\in\RR^n:\bm{Ax}\leq \bm{b}\}$,
	we can generate samples on its interior using the following procedure:
	\begin{enumerate}[label=(\roman*)]
		\item Pick a starting point $\bm{x}_0\in\RR^n$
		which lies in the interior of the polyhedron
		$\mathcal{P}$.
		\item Generate a unit random direction
		$\bm{u}\in\RR^n$
		\item Calculate the minimum and maximum
		values of $\lambda\in\RR$ such that
		$\bm{x}_0+\lambda\bm{u}\in\mathcal{P}$.
		Those values $\lambda_{\min}$ and $\lambda_{\max}$ can be calculated in close
		form by requring that $\lambda \bm{u}\leq \bm{b}-\bm{Ax}$.
		Note that for general polyhedras, it can happen that $\lambda_{\max}=+\infty$ or $\lambda_{\min}=-\infty$,
		but the polyhedras we examine are all bounded,
		so $\lambda_{\max}$ and $\lambda_{\min}$ are finite.
		\item Pick a $\lambda_*$
		uniformly at random from the set $[\lambda_{\min},\lambda_{\max}]$ and
		generate sample $\bm{x}_*=\bm{x}_0+\lambda_*\bm{u}$.
		\item Repeat the procedure from Step (ii) using 
		$\bm{x}_*$ as the new starting point.
		Terminate whenever we generate the required number of samples.
	\end{enumerate}
	We then follow this procedure to generate
	samples for all polyhedra $\mathcal{P}\in\mathcal{C}$.
	If the set $\mathcal{P}\in\mathcal{C}$ contains
	$2$ polyhedra with the same representation more than
	once, then we only sample one of them.\par 
	The purpose of this procedure is to generate
	samples in areas where there is classification 
	ambiguity by the tree learners.
	Through those samples, we attempt to improve the generalization
	of ML models
	and thus create better constraint approximators.
	This sampling procedure is adaptive and
	can be repeated to iteratively generate more refined samples
	and more accurate approximations.
\end{enumerate}

\subsection{Enhancement 3: Relaxations}
As we described in previous sections, in order
to solve global Optimization Problems, we first
create a MIO approximation of the original
problem and we then optimize over
this MIO in order to receive an approximate
solution $\bm{x}_*$.
Then, we perform a local-search PGD step
in order to improve the solution $\bm{x}_*$
in terms of both feasibility and optimality.\par 
However, due to the nature of this
approach, it may occur that the MIO approximation
is not feasible, even if the original
problem is feasible.
For instance, if a constraint is very difficult to
satisfy, then we may have very limited feasible samples of
that constraint.
Thus, the learner may approximate
the constraint as infeasible almost everywhere,
leading to an infeasible MIO.

%
%
To mitigate this issue, we relax the approximating MIO so
that it will always have a feasible solution, provided that the
original problem is also feasible.
In particular, given the approximation MIO of Eq. (\ref{eq:mio}),
we relax it as follows:
\begin{equation}\label{eq:mio_rel}
	\begin{aligned}
		\min_{\bm{x}\in\RR^n,u\in\RR} \quad & y_{ML\_REGR}(\bm{x})+\lambda \sum_{i\in I}u_i\\
		\textrm{s.t.} \quad &y^{(i)}_{ML}(\bm{x})+u_i\geq a_i,\: i\in I,\\
		&\bm{Cx}\leq  \bm{d},\\
		&u_i\geq 0,\quad i\in I,
	\end{aligned}
\end{equation}
where $u_i$ are the relaxation variables and $\lambda$ is a parameter that 
determines how much to penalize relaxed constraints.
By including this relaxation variable in the constraints,
we are ensuring that the MIO approximate constraints do not introduce any 
infeasibilities.
Also, the more we increase the penalty $\lambda$,
the more the MIO solver will prioritize 
solutions with small values of the variable $u_i$ and thus less relaxed
constraints will be.
However, we also have to note that we first attempt to solve the problem
without relaxations, and if the approximating MIO is infeasible,
then we resolve using the described relaxation technique.

\subsection{Enhancement 4: Robustness}
\label{section:ro}
When using ML models to approximate nonlinear functions through
samples, there can be a great deal of uncertainty in the learned model
parameters.
In particular, training on a different set of samples
can lead to different model parameters and thus different 
approximations of the nonlinear function. 
To partially account for this issue, we 
will use Robust Optimization when embedding
the ML learners into the final MIO.\par 
Robust Optimization (RO) is a methodology for dealing 
with uncertainty in the data of an optimization problem 
\cite{ben-tal_robust_2009}.
In our approach, however, we will use RO
to deal with uncertainty not in the model data,
but in the trained model parameters.
For all types of models, we chose to model
uncertainty using the $p$-norm uncertainty set
$\mathcal{U}^{\rho}_p=\{\bm{x}\in\RR^n:||\bm{x}||_p\leq \rho\}$.
The precise way of defining the uncertainty is 
slightly different across the learners, so
we will describe each of the different models separately:
\begin{itemize}
	\item \textit{Support Vector Machines}:
	In the case of SVMs, we attempt to model the fact
	that after training the SVM, there is uncertainty in the
	values of the coefficient vector $\bm{\beta}$.
	For our use-case, we will choose to model $\bm{\beta}$ using
	multiplicative uncertainty,
	in order to scale the coefficients of $\bm{\beta}$ proportionally
	to their nominal values.
	More concretely, 
	if an SVM is used to approximate a constraint $g(\bm{x})\leq 0$
	and $\bar{\beta}_0$ and $\bar{\bm{\beta}}$ are the nominal parameters of the SVM
	after training, 
	then we approximate $g(\bm{x})\leq 0$ with
	the uncertain constraint:
	\begin{equation}
		\bar{\beta}_0+(\bar{\bm{\beta}}\odot(\bm{1}+\bm{z}))^T\bm{x}\geq 0,\quad \forall \bm{z}\in\mathcal{U}^{\rho}_p,
	\end{equation}
	where $\odot$ denotes the element-wise vector product and $\bm{1}$ is the vector of ones.
	Then, by using the tools presented in \cite{bertsimas_robust_2022}, this constraint can be written equivalently
	as:
	\begin{equation}
		\bar{\beta}_0+\bar{\bm{\beta}}^T\bm{x}-\rho ||\bar{\bm{\beta}}\odot \bm{x}||_q\geq 0,
	\end{equation}
	where $||\cdot ||_q$ is the dual norm of $||\cdot ||_p$ (i.e., $\frac{1}{q}+\frac{1}{p}=1$).
	
	\item \textit{Decision Trees - Gradient Boosted Machines}:
	In the case of Decision Trees,
	we will model the uncertainty in the coefficient vectors $\bm{a}_j$ of the hyperplane
	splits.
	As with SVMs, we will again use multiplicative uncertainty to scale the
	coefficients proportionately to their nominal values.
	In this case, multiplicative uncertainty is crucial,
	as it allows us to keep the zero elements of the coefficient vectors constant.
	For instance, if the
	original Decision Tree uses parallel splits, 
	we only want to consider uncertain vectors $\bm{a}_j$ that
	represent parallel splits, an effect which will be captured
	by using multiplicative uncertainty.\par 
	More formally, if we want to approximate the constraint $g(\bm{x})\leq 0$
	with a Decision Tree and $\bar{\bm{a}}_j$ is the coefficient vector of the $j$
	node of the tree after training, then the splitting constraints of Eq. (\ref{eq:dt})
	will take the following form after accounting for uncertainty:
	\begin{equation}
		\begin{aligned}
			(\bar{\bm{a}}_j\odot (\bm{1}+\bm{u}))^T\bm{x}&\leq b_j + M(1-z_i),\quad \forall \bm{u}\in\mathcal{U}_p^{\rho}, \quad \forall i\in\mathcal{L},j\in L(L_i),\\
			(\bar{\bm{a}}_j\odot (\bm{1}+\bm{u}))^T\bm{x}&\geq  b_j - M(1-z_i)+\epsilon, \quad \forall \bm{u}\in\mathcal{U}_p^{\rho},\quad \forall i\in\mathcal{L},j\in R(L_i),
		\end{aligned}
	\end{equation} 
	Again, by using the tools in \cite{bertsimas_robust_2022}, we can write the above constraints
	equivalently as:
	\begin{equation}\label{eq:dt_ro}
		\begin{aligned}
			\bar{\bm{a}}_j^T\bm{x}+\rho ||\bar{\bm{a}}_j\odot \bm{x}||_q&\leq b_j + M(1-z_i), \quad \forall i\in\mathcal{L},j\in L(L_i),\\
			\bar{\bm{a}}_j^T\bm{x}-\rho ||\bar{\bm{a}}_j\odot \bm{x}||_q&\geq  b_j - M(1-z_i)+\epsilon,\quad \forall i\in\mathcal{L},j\in R(L_i),
		\end{aligned}
	\end{equation} 
	where $||\cdot ||_q$ is the dual norm of $||\cdot ||_p$.\par 
	Finally, in the case of Gradient Boosted
	Machines, 
	we use robustness by enforcing the constraints
	of Eq. (\ref{eq:dt_ro}) to each one of the
	individual tree learners of the GBM ensemble.
	We make sure to not use too many learners in our GBM 
	ensemble, in order to avoid infeasibilities that
	arise from enforcing robustness in all learners at the same time. 
\end{itemize}
\bigskip\bigskip
\section{An Illustrative example}
To better illustrate the different steps of
the enhanced optimization methodology, we consider the 
following non-convex problem:
\begin{equation}\label{eq:example}
	\begin{aligned}
		\min_{\bm{x}\in\RR^n} \quad &f(\bm{x})=-x_2\\
		\textrm{s.t.} \quad &g_1(\bm{x})=-0.43\ln(x_1-0.5)-1.1-x_1+x_2\leq 0,\\
		&g_2(\bm{x})=-x_2+0.33\ln(x_1-0.4)+1.2-0.2x_1\leq 0,\\
		&g_3(\bm{x})=-x_2+1.1x_1+0.3\leq 0,\\
		&g_4(\bm{x})=-x_2-1.5x_1+2.6\leq 0,\\
		&0.51\leq x_1\leq 1.5,\\
		&0.3\leq x_2\leq 1.6.
	\end{aligned}
\end{equation}
To facilitate understanding of the method,
we chose a problem with two variables and a linear 
objective, although the methodology also supports non-linear
and non-convex objectives.

\subsection{Standard Form Generation}

The first step is to separate the linear from
the non-linear constraints and bring the problem into the standard 
form described in Eq. (\ref{eq:gop_std}).
After doing that, Problem (\ref{eq:example}) becomes:
\begin{alignat*}{3}
	\min& \:f(\bm{x})=-x_2 & \text { Objective } \\
	\hline \text { s.t. } & g_1(\bm{x})=-0.43\ln(x_1-0.5)-1.1-x_1+x_2\leq 0, & \text { Nonlinear } \\
	&g_2(\bm{x})=-x_2+0.33\ln(x_1-0.4)+1.2-0.2x_1\leq 0, & \text { constraints } \\
	\hline
	&g_3(\bm{x})=-x_2+1.1x_1+0.3\leq 0, & \text { Linear } \\
	&g_4(\bm{x})=-x_2-1.5x_1+2.6\leq 0, & \text { constraints } \\
	\hline 
	&0.51\leq x_1\leq 1.5, & \text { Variables } \\
	&0.3\leq x_2\leq 1.6. & \text { and bounds }
\end{alignat*}
\noindent
The linear constraints and the linear objective
are directly passed to the MIO model.
The nonlinear constraints involving $g_1$ and $g_2$ will be
approximated in the next steps of the methodology.
However, before proceeding into the next steps, we make
sure that all variables of the nonlinear constraints
are bounded.
If they are not, then we attempt to compute bounds the way that is described in
Section \ref{section:octh}.

\subsection{Sampling of Nonlinear constraints}
In this step, we sample the nonlinear constraints
$g_1(\bm{x})\leq 0$ and $g_2(\bm{x})\leq 0$.
The goal is to obtain samples $\{(\bm{\widetilde x_k},\mathds{1}\{g_1(\bm{\widetilde x_k})\leq 0\})\}^N_{k=1}$
and $\{(\bm{\widetilde x_k},\mathds{1}\{g_2(\bm{\widetilde x_k})\leq 0\})\}^N_{k=1}$
that will be used to train ML models for approximating the nonlinear constraints $g_1(\bm{x})\leq 0$
and $g_2(\bm{x})\leq 0$.\par 
The sampling is performed in the following steps: (i) Boundary Sampling,
(ii) Latin Hypercube Sampling, (iii) kNN quasi-newton sampling
and (iv) OCT-based adaptive sampling. 
The first $3$ sampling steps are the same as the ones used in OCTHaGOn and
are described in detail in Section \ref{section:octh}.
Then, Sampling Step (iv), which is part of our enhancements, is used to adaptively
obtain more refined samples in areas where the nonlinear constraints
are not approximated well by the ML learners (e.g. near the feasibility
boundaries of the constraint).
This step is described in detail in section \ref{section:gop_sampling}.\par
\begin{figure}[h]
	\centering
	\begin{subfigure}[b]{0.49\textwidth}
		\centering
		\includegraphics[width=\textwidth]{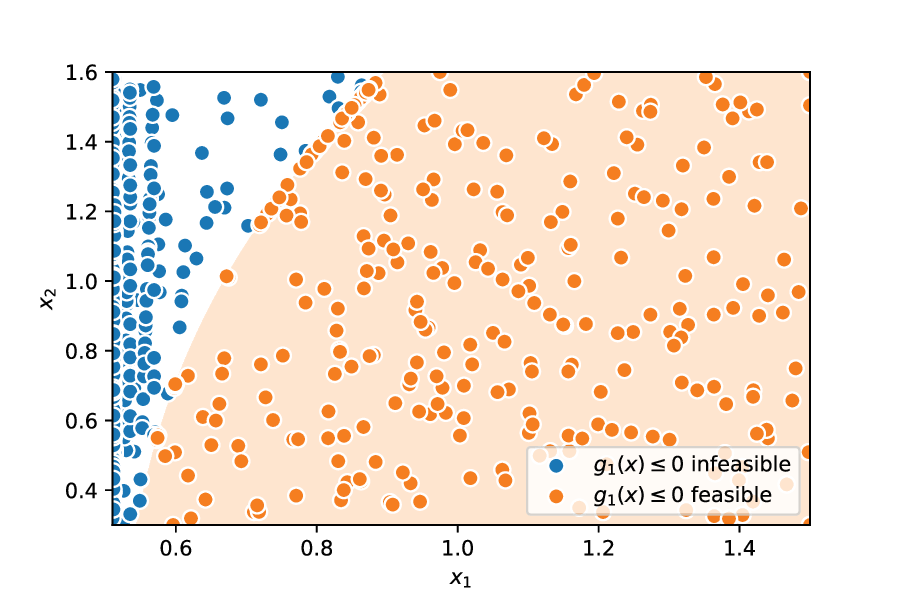}
		\caption{Samples of $g_1(\bm{x})\leq 0$.}
		\label{fig:g1}
	\end{subfigure}
	\hfill
	\begin{subfigure}[b]{0.49\textwidth}
		\centering
		\includegraphics[width=\textwidth]{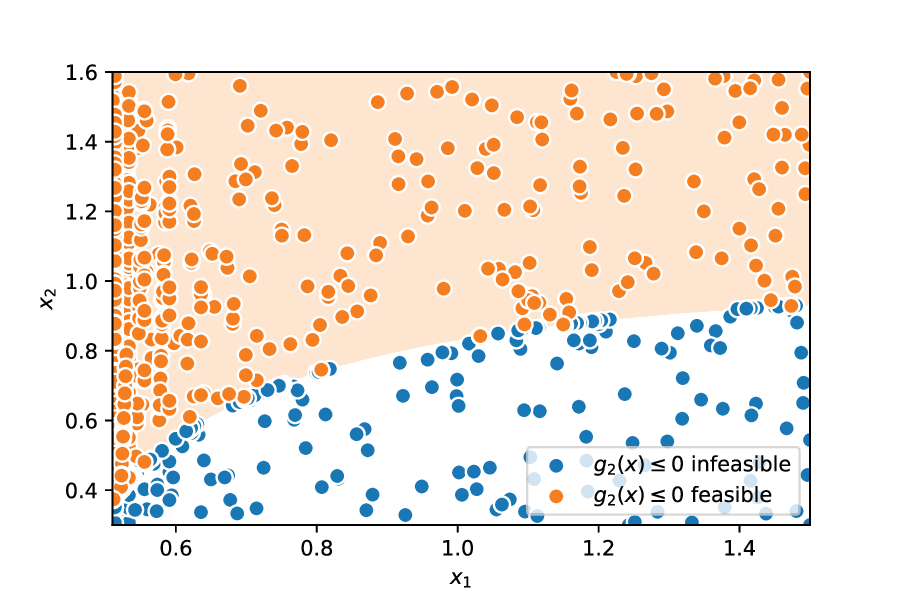}
		\caption{Samples of $g_2(\bm{x})\leq 0$.}
		\label{fig:g2}
	\end{subfigure}
	\caption{Feasible region and samples of the nonlinear constraints.}
	\label{fig:g1g2}
\end{figure}

 
After applying those sampling steps, we can see the samples
obtained for the constraints $g_1(\bm{x})\leq 0$ and $g_2(\bm{x})\leq 0$ in
the Figures \ref{fig:g1} and \ref{fig:g2},
respectively.
By examining the samples and the orange feasibility regions
of the constraints, we can see that the adaptive sampling procedures
we used has helped generate samples near the boundaries of the 
constraint feasibility regions, which is essential for
good ML approximations.

\subsection{Model Training}
Given the samples we obtained from the 
previous step,
we train a number of learners  in order to approximate the nonlinear constraints $g_1(\bm{x})\leq0$ and
$g_2(\bm{x})\leq 0$.
The learners are trained on the datasets
$D_1=\{(\bm{\widetilde x_k},\mathds{1}\{g_1(\bm{\widetilde x_k})\leq 0\})\}^N_{k=1}$
and 
$D_2=\{(\bm{\widetilde x_k},\mathds{1}\{g_2(\bm{\widetilde x_k})\leq 0\})\}^N_{k=1}$
that contain the samples for the constraints $g_1(\bm{x})\leq 0$ and $g_2(\bm{x})\leq0$ respectively.
For each one of those datasets, we separated the samples into a training and a validation set.
We then used the training sets to train a Multi-layer perceptron (MLP), a Support Vector Machine (SVM),
a Gradient Boosted Machine (GBM) and Hyperplane-based decision tree (OCT).
Then, we measured the accuracy of those models in the respective validation sets
and for each of the $2$ constraints, we picked the learner that demonstrated
the highest accuracy.\par
\begin{figure}[h]
	\centering
	\begin{subfigure}[b]{0.49\textwidth}
		\centering
		\includegraphics[width=\textwidth]{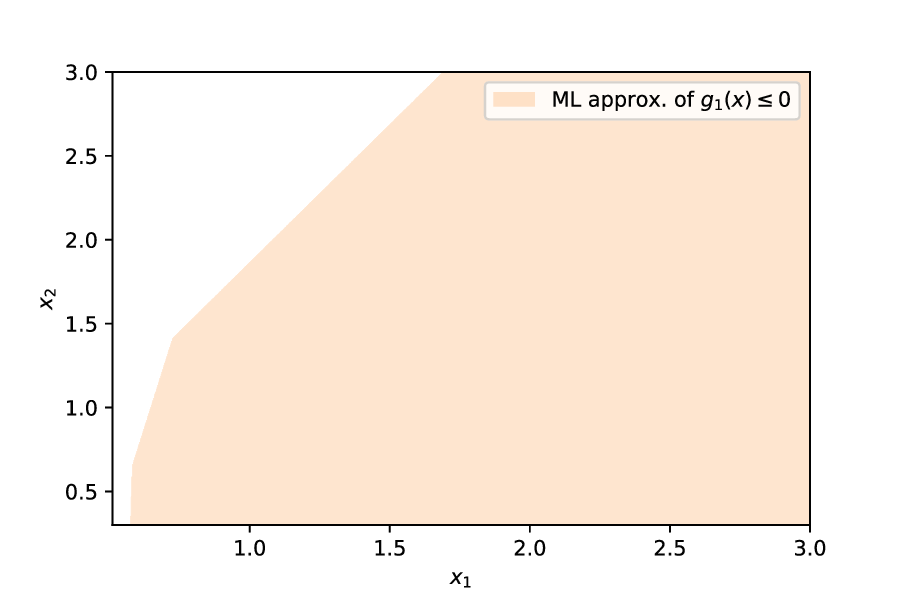}
		\caption{MLP approximation of $g_1(\bm{x})\leq 0$.}
		\label{fig:g1_approx}
	\end{subfigure}
	\hfill
	\begin{subfigure}[b]{0.49\textwidth}
		\centering
		\includegraphics[width=\textwidth]{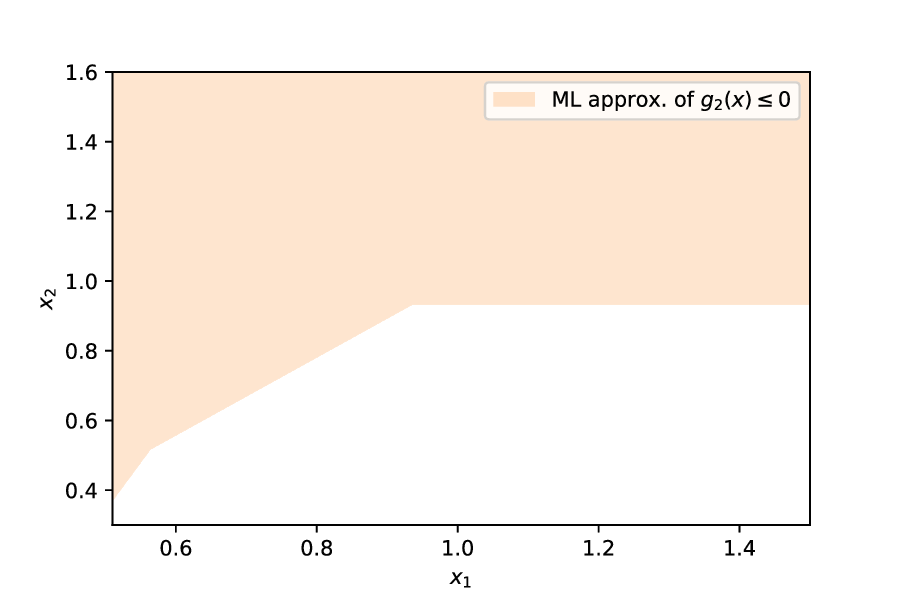}
		\caption{OCT approximation of $g_2(\bm{x})\leq 0$.}
		\label{fig:g2_approx}
	\end{subfigure}
	\caption{ML approximation of the nonlinear constraints.}
	\label{fig:g1g2_approx}
\end{figure}

 
In the case of the constraint $g_1(\bm{x})\leq 0$, the best learner was the MLP
with an accuracy of $0.97$, whereas in the case of the constraint $g_2(\bm{x})\leq 0$,
the best learner was the OCT with an accuracy of $0.99$.
The ML approximation for the constraints involving $g_1$
and $g_2$ are shown in Figures \ref{fig:g1_approx} and \ref{fig:g2_approx},
respectively.

\subsection{MIO representation}
\begin{figure}[h]
	\begin{minipage}{.45\textwidth}
		\includegraphics[width=0.9\textwidth]{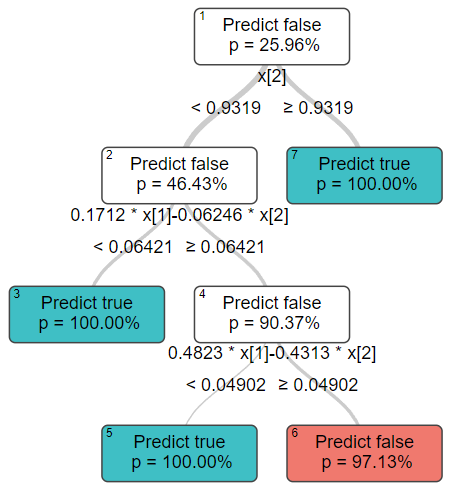}
	\end{minipage}%
	\begin{minipage}{.45\textwidth}
		\begin{equation*}\label{eq:oct_example}
			\begin{aligned}
				y_{OCT} &= z_1+z_2+z_3,\\
				1&=z_1+z_2+z_3+z_4,\\
				x_2&\geq 0.9319-M(1-z_1),\\
				x_2&\leq 0.9319+M(1-z_2)-\epsilon,\\
				0.1712x_1-0.06246x_2&\leq 0.06421+M(1-z_2)-\epsilon,\\
				x_2&\leq 0.9319+M(1-z_3)-\epsilon,\\
				0.1712x_1-0.06246x_2&\geq 0.06421-M(1-z_3),\\
				0.4823x_1-0.4313x_2&\leq 0.04902+M(1-z_3)-\epsilon,\\
				x_2&\leq 0.9319+M(1-z_4)-\epsilon,\\
				0.1712x_1-0.06246x_2&\geq 0.06421-M(1-z_4),\\
				0.4823x_1-0.4313x_2&\geq 0.04902+M(1-z_4),\\
			\end{aligned}
		\end{equation*}
	\end{minipage}
	\caption{OCT approximator of $g_2(\bm{x})\leq 0.$}
	\label{fig:ml_appro_g2}
\end{figure}

In this step, we represent the trained learners
using an MIO formulation as described in Section \ref{section:ml-models}.
In particular, the OCT approximator of the constraint $g_2(\bm{x})\leq 0$
can be shown in Figure \ref{fig:ml_appro_g2} in the form of a hyperplane-based
decision tree.
The tree consists of $3$ non-terminal and $4$ terminal (leaf) nodes.
The leaf nodes that correspond to feasible regions
are the ones shown in blue.\par 
In order to represent this decision tree
using an MIO formulation, we introduce $4$ binary auxiliary variables
($z_1,z_2,z_3,z_4$), 1 for each leaf.
Each auxiliary variable becomes active if and only if $\bm{x}$ lies in the corresponding
leaf.
Then, we use a big-M formulation to encode the output of the
Decision Tree, as shown next to Figure \ref{fig:ml_appro_g2}.
Note that in the resulting formulation, $y_{OCT}$ is the output
of the decision tree and $\epsilon$ is a small positive
constant used to model strict inequalities.
Then, given the MIO representation of the decision tree of Figure \ref{fig:ml_appro_g2},
we can approximate the constraint $g_2(\bm{x})\leq 0$ with the constraint $y_{OCT}\geq 0.5$.\par 

\begin{figure}[h]
	\begin{minipage}{.5\textwidth}
		\includegraphics[width=0.8\textwidth]{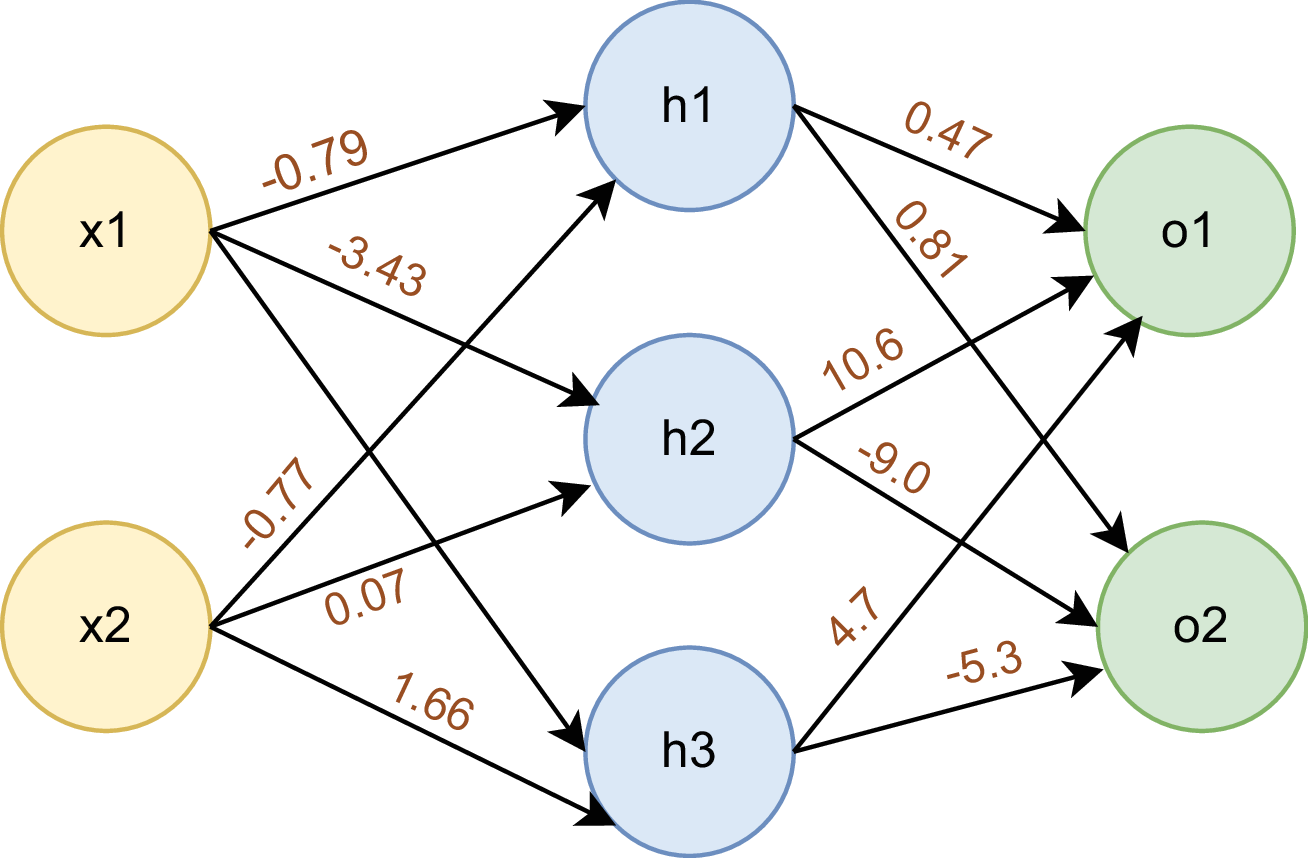}
	\end{minipage}%
	\begin{minipage}{.5\textwidth}
		\begin{equation}\label{eq:mlp_example}
			\begin{aligned}
				y_{MLP}&=\bm{a_3}^T\cdot [-1,\:1], \\
				\bm{a}_3&\geq \bm{0},\\
				\bm{a}_3&\geq \bm{W_2}\bm{a_2}+\bm{b_1},\\
				\bm{a}_3&\leq M \mathbf{v_2},\\
				\bm{a}_3&\leq \bm{W_2}\bm{a_2}+\bm{b_2}-M(1-\bm{v_2}),\\
				\bm{a}_2&= \bm{W_1}\bm{x}+\bm{b_1},\\
				\bm{a}_2&\in\RR^3,\:\bm{a}_2\in\RR^3,\: \bm{v_2}\in\{0,1\}^3.
			\end{aligned}
		\end{equation}
	\end{minipage}
	\caption{MLP approximator of $g_1(\bm{x})\leq 0.$}
	\label{fig:ml_appro_g1}
\end{figure}

Next, the trained MLP that approximates the constraint $g_1(\bm{x})\leq 0$
is shown in Figure \ref{fig:ml_appro_g1}.
This MLP consists of $2$ input nodes,
$3$ hidden nodes with ReLU activations and $2$ output
nodes.
It was trained using Negative Log Likelihood loss
in a binary classification task, so that it predicts that
the constraint is feasible whenever $o_2\geq o_1$ (i.e., $o_1$ and $o_2$
are the output nodes shown in Figure \ref{fig:ml_appro_g1}).
Note here that although we have a binary classification task, we have used
$2$ output nodes instead of $1$, where an output pair $(o_1,o_2)=(1,0)$ 
corresponds to an infeasible input $\bm{x}$, while an output pair
$(o_1,o_2)=(0,1)$ corresponds to a feasible input $\bm{x}$. \par 
In order to model the MLP using MIO, we first use a vector of
continuous variables $\bm{a}_i$ to represent the input of
the $i$-th layer of the MLP (i.e., the input of the first hidden layer is denoted as $\bm{a}_2$,
the input of the next layer as $\bm{a}_3$ etc.).
Then, for the hidden layer, we additionally use a vector of binary
variables $\bm{v_2}$ to model the ReLU activation.
The resulting big-M formulation is described in Eq. (\ref{eq:mlp_example}),
where $\bm{W_1},\bm{W_2}$ are the weight matrices 
and $\bm{b_1},\bm{b_2}$ are the bias matrices of the trained MLP.
Using this notation, the constraint $g_1(\bm{x})\leq 0$ can be approximated
with the constraint $y_{MLP}\geq 0$, where $y_{MLP}$ is given in Eq. (\ref{eq:mlp_example}).

\begin{figure}[H]
	\centering
	\includegraphics[width=0.7\textwidth]{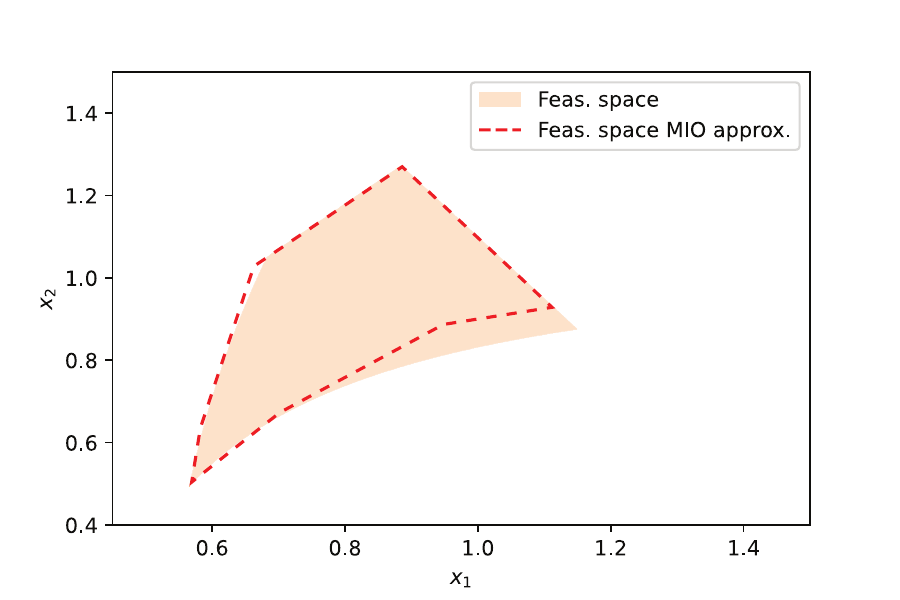}
	\caption{MIO approximation of feasible space.}
	\label{fig:feas_approx}
\end{figure}

Then, we combine the MIO formulations of the OCT and the MLP into
a unified MIO formulation. 
In this formulation, we also include the linear constraints of the original problem,
and this way we end up with an MIO approximation of the problem.
In Figure \ref{fig:feas_approx} we can see the resulting MIO
approximation of the feasible space against the actual feasible
space of the original problem.
We can see that the non-convexity of the feasible region is captured
by the MIO formulation, and thus by optimizing over this MIO
approximation, we are (approximately) solving the original problem.

\subsection{Robustness} 
Constructing a data-driven MIO aproximation of the original
problem is not always accurate, in part due to the uncertainty
introduced by the sampling process.
This uncertainty be partially accounted for by introducing
a level of robustness into the MIO formulation.
For example, following the steps we described in Section \ref{section:ro}, 
we can rewrite the MIO representation of the OCT-H of Figure \ref{fig:ml_appro_g2} in the following robust way:
\begin{equation}\label{eq:oct_example_ro}
	\begin{aligned}
		y_{OCT} &= z_1+z_2+z_3,\\
		1&=z_1+z_2+z_3+z_4,\\
		x_2-{\rho ||x_2||_q}&\geq 0.9319-M(1-z_1),\\
		x_2+{\rho ||x_2||_q}&\leq 0.9319+M(1-z_2)-\epsilon,\\
		0.1712x_1-0.06246x_2+\rho||0.1712x_1-0.06246x_2||_q&\leq 0.06421+M(1-z_2)-\epsilon,\\
		x_2+\rho ||x_2||_q&\leq 0.9319+M(1-z_3)-\epsilon,\\
		0.1712x_1-0.06246x_2-\rho||0.1712x_1-0.06246x_2||_q&\geq 0.06421-M(1-z_3),\\
		0.4823x_1-0.4313x_2+\rho||0.4823x_1-0.4313x_2||_q&\leq 0.04902+M(1-z_3)-\epsilon,\\
		x_2+\rho ||x_2||_q&\leq 0.9319+M(1-z_4)-\epsilon,\\
		0.1712x_1-0.06246x_2-\rho||0.1712x_1-0.06246x_2||_q&\geq 0.06421-M(1-z_4),\\
		0.4823x_1-0.4313x_2-\rho||0.4823x_1-0.4313x_2||_q&\geq 0.04902+M(1-z_4).\\
	\end{aligned}
\end{equation}
Next, depending on the value of $q$ (which is determined
by the type of uncertainty we are protecting against), we 
reformulate the norm operators in a tractable way.
For instance, if we want to protect against $L_1$ or $L_{\infty}$ uncertainty
sets, then $q$ takes the value $\infty$ and $1$, respectively,
and the robust MIO can be easily reformulated using linear constraints.
On the other hand, if we want to protect against $L_2$ uncertainty,
then the resulting MIO is conic quadratic.\par 
Note that the extend to which we include robustness
into our MIO approximation is determined by the hyperparameter $\rho$:
A value of $\rho=0$ indicates that we do not want to
consider robustness, while higher values of $\rho$ correspond
to bigger uncertainty sets and more conservative robust counterparts.
For the purpose of this example, we use a value of $\rho=0.1$.
\subsection{Solve MIO \& Improve} 
After having constructed an MIO approximation of Problem \ref{eq:example}, 
we optimize over this approximation using Gurobi.
The initial solution we get is $[x_1,x_2]=[1.108,0.937]$ with
an objective value of $-1.108$.
This initial solution is not optimal, 
due to the MIO approximation error shown in
Figure \ref{fig:feas_approx}.\par
\begin{figure}[H]
	\centering
	\includegraphics[width=0.7\textwidth]{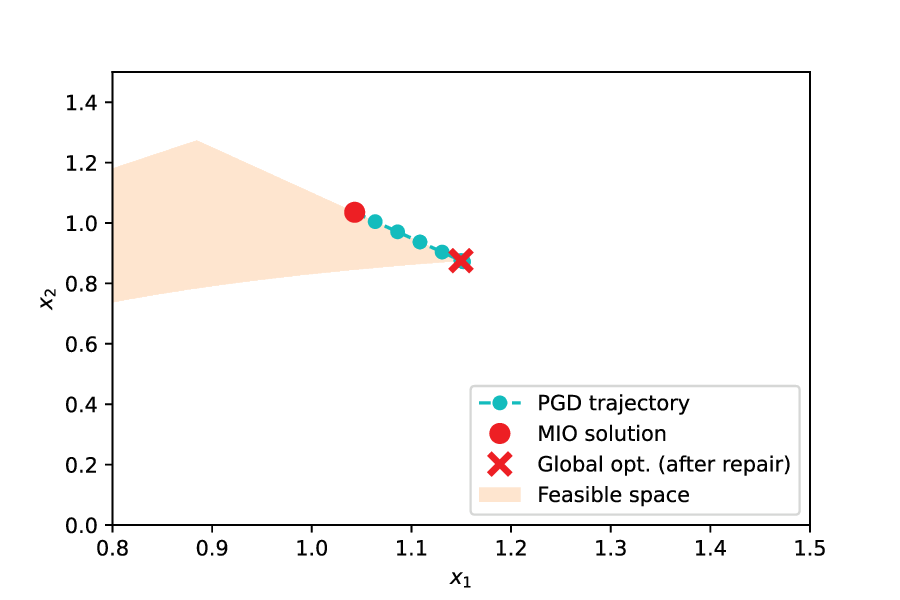}
	\caption{PGD trajectory.}
	\label{fig:pgd}
\end{figure}

 
We then proceed to improve the solution using
the PGD local-search procedure described in Section \ref{section:octh}.
Note that in this procedure, we also introduce 
momentum and we conditionally utilize second order information
for faster and more accurate convergence.
The resulting solution after $10$ iterations of the solution improvement procedure
is $[x_1,x_2]=[1.1497,0.875]$ with an objective value of $-1.1497$,
which is the global optimum.
The trajectory followed by the PGD iterations can be seen
visually in Figure \ref{fig:pgd}.

\section{Computational Experiments On Benchmarks}
\label{section:experiments}
In order to test the new framework, which we will refer to as GoML, we use a number
of global optimization problems from the literature.
We focus on continuous non-convex optimization problems
where the optimization variables are all bounded.
The reason we concentrate on problems with bounded
variables is that GoML and OCTHaGOn both
require variable bounds in order to perform the
initial sampling steps.
In case variable bounds are not provided, the frameworks
try to infer those bounds, but
we will only focus on already bounded problems
in order to have a more consistent ground of comparison.

\subsection{Experimental Procedure}
In order to perform the experiments, we use two different versions of
our GoML framework.
In the first version (GoML Base),
we use the standard convex version of Gurobi
to solve the MIO approximation of the original
problem.
This means that we pass directly to Gurobi
the linear constraints (and objectives) of the original problem,
and we approximate all the non-linear ones using ML models.
This approach is analogous to the one used by OCTHaGOn.
In the second version (GoML NonCVX),
we also pass the convex constraints and the non-convex quadratic constraints
straight to Gurobi, and we use ML to approximate only the non-convex constraints
that are not quadratic.
The second approach leverages the ability of Gurobi to
solve non-convex MINLP problems, as long as the non-convex constraints
are quadratic.
The caveat of this approach is that 
some of the benchmarks contain only non-convex
quadratic constraints and objectives, and thus in those instances,
GoML NonCVX offloads all the constraints (and objective) to Gurobi.
However, in many of those instances, GoML Base is also
able to solve them equally well, as seen on the result section.\par 
Both versions of GoML have a number of hyperparameters, the most important of
which is the robustness radius $\rho$ and the relaxation coefficient $\lambda$.
In order to determine the best value of those hyperparameters, we perform
a grid search.
In particular, we run GoML for different values of the hyperparameters
and we keep the solution with the best objective value.
For the robustness parameter $\rho$, we use
values $\rho=0$ (i.e., no robustness)
and $\rho=0.01,0.1,1$ which represent uncertainty sets
of different sizes.
For the relaxation coefficient $\lambda$, we first
experiment with no relaxations (corresponding to $\lambda=\infty$)
and then, with varying levels of relaxation penalties
($\lambda=10^{2}$ and $\lambda=10^4$).\par 
If we use a naive grid search approach and 
we solve the problem from scratch for each
combination of hyperparameters, then our approach will
perform very poorly time-wise.
However, we notice that when we change the robustness 
and relaxation parameters, we do not  need to retrain the
learners, since both robustness and relaxations
are applied to the MIO approximations after the
learners are trained.
Hence, in our grid search, we only need to sample the constraints
and train the ML models once, and then
we can re-embed the models and resolve the MIO
for each set of hyperparameters.
This observation is really important, since the sampling
and training steps take the vast majority of the solution
time (i.e., in many cases, more than $95\%$).\par

\subsection{Results}
Using grid-search the way we described, we use GoML Base
and GoML NonCVX to solve a
range of continuous global optimization problems
from MINLPLib \cite{bussieck_minlplibcollection_2003-1}.
For our benchmarks, we use problems that contain
between $1$ and $110$ variables and between $1$ and $88$ constraints.
We first benchmark our methodologies (GoML base and GoML NonCVX) against OCTHaGOn,
in order to see the effect of our enhancements.
Then, we compare both approaches against BARON \cite{sahinidis_baron_1996},
a well-established global optimizer for mixed-integer nonlinear problems.
For all methods, we measure the percent optimality gap and
the solution time in seconds.\par 
We set a pre-specified time limit of $1500$s for all solvers,
and we run our experiments in a Dell Inspiron Laptop
with an 8-core Ryzen 5700U processor clocked at 1.8GHz
and 16GB of RAM.
The enhanced GoML framework is implemented in Julia 1.6.2,
and is an extension of the publicly available OCTHaGOn 
code-base \cite{bertsimas_global_2022}.
In our experiments, we use Gurobi 8.1 as the underlying MIO solver for both GoML and OCTHaGOn, 
and we compare against BARON version 2021.01.13.

\par 
The comparison results are shown in Table \ref{tab:results}.
In this table, we have recorded the name of the problem,
the number of variables, the optimality gap and timing
for the different methods.
When the optimality gap is below $0.1\%$,
we consider the solution optimal and we record ``GOpt"
in the table.
In the case of GoML NonCVX, we use an asterisk $*$
to label the instances where the problems are non-convex quadratics
and thus all the constraints and the objective are directly passed
into Gurobi.
Also, when a method fails to produce
a feasible solution to the problem,
we record "Inf" in the respective
columns.\par 
By looking at the results of Table \ref{tab:results},
we can see that although GoML is an approximate
approach for solving global optimization problems,
it is able to solve many of the problems to near-optimality.
Also, both GoML Base and GoML NonCVX are at least as good as 
OCTHaGOn in the vast majority of instances.
In particular,
by comparing GoML Base with OCTHaGOn, we can see that
GoML Base has better optimality gaps in $39$ out of the $77$ instances,
and worse optimality gaps in only $4$ of the $77$ instances.
This means that GoML Base offers solutions that are at least as good 
as OCTHaGOn in $73$ out of the $77$ instances, validating that our enhancements
are indeed effective.\par 

On the other hand, BARON is able to solve $73$ out of the $77$ instances
to optimality, which is expected since those problems have long been used
as benchmarks in the global optimization literature.
However, in $4$ of the $77$ instances, GoML Base provides solutions with 
better optimality gaps than BARON, whereas in another $3$ of the $77$ instances,
both GoML Base and BARON solve the problem to optimality, but GoML base
has better solution times.



Alternatively, if we use GoML NonCVX instead of GoML Base in our comparisons,
then we will see even better results in terms of both optimality gap
and solution time.
In particular, by offloading compatible nonlinear constraints to Gurobi,
we improve the optimality gaps of GoML Base in $21$ of the $77$ instances,
and we also improve solution time in $50$ out of the $51$ instances
solved to optimality by both GoML Base and GoML NonCVX.
However, the caveat of GoML NonCVX is that there
are instances where GoML NonCVX finds the optimal solution
by offloading all the constraints to Gurobi (i.e.,
those instances are labeled with an asterisk $*$), effectively not utilizing the other
parts of the methodology.
Nevertheless, GoML Base also seems to be able to find the 
optimal solution in many of those instances.\par 
Finally, if we take the best solutions of GoML Base and GoML NonCVX,
while excluding solutions labeled with an asterisk (*),
then there are $9$ instances where the enhanced GoML framework
has better optimality gap or solution time than BARON.
In particular, in $4$ of those instances, the enhanced framework
offers better optimality gaps than BARON, and in another $5$
instances, both BARON and GoML solve the problem to optimality,
but GoML has a better solution time.
The corresponding results are all shown in Table \ref{tab:results}.

\clearpage
\captionsetup[longtable]{labelfont=bf,textfont=bf,position=bottom,font=small}
{\renewcommand{\arraystretch}{1.08}
	\setlength{\tabcolsep}{3pt}		
		\begin{longtable}{|c|c|cc|cc|cc|cc|}
			
			\hline 
			&&\multicolumn{2}{c|}{GoML Base} &\multicolumn{2}{c|}{GoML NonCVX} & \multicolumn{2}{c|}{BARON} & \multicolumn{2}{c|}{OCTHaGOn}\\
			name & \# vars & gap (\%) &  time (s) & gap (\%) &  time (s) & gap (\%) &  time (s) & gap (\%) &  time (s) \\
			\hline \hline
			\endfirsthead
			\hline 
			&&\multicolumn{2}{c|}{GoML Base} &\multicolumn{2}{c|}{GoML NonCVX} & \multicolumn{2}{c|}{BARON} & \multicolumn{2}{c|}{OCTHaGOn}\\
			name & \# vars & gap (\%) &  time (s) & gap (\%) &  time (s) & gap (\%) &  time (s) & gap (\%) &  time (s) \\
			\hline \hline
			\endhead
			ex4\_1\_7 &      1 & \textbf{GOpt} &            $7.4$ &          $0.12$ &               $7$ & \textbf{GOpt} &  $\mathbf{0.1}$ &        $16.5$ &            $1.2$ \\
			ex4\_1\_2 &      1 & \textbf{GOpt} &            $7.3$ &   \textbf{GOpt} &               $6$ & \textbf{GOpt} &  $\mathbf{0.4}$ &        $4.58$ &              $1$ \\
			ex4\_1\_3 &      1 & \textbf{GOpt} &            $8.1$ &   \textbf{GOpt} &             $5.9$ & \textbf{GOpt} &  $\mathbf{0.1}$ &        $2.31$ &            $1.3$ \\
			ex4\_1\_6 &      1 & \textbf{GOpt} &            $6.9$ &   \textbf{GOpt} &             $5.1$ & \textbf{GOpt} &  $\mathbf{0.2}$ &      $607.74$ &            $1.3$ \\
			ex4\_1\_1 &      1 & \textbf{GOpt} &            $119$ &   \textbf{GOpt} &           $116.5$ & \textbf{GOpt} &  $\mathbf{1.2}$ &         $143$ &           $11.7$ \\
			prob06 &      1 & \textbf{GOpt} &           $71.1$ &  \textbf{GOpt*} &             $1.7$ & \textbf{GOpt} &  $\mathbf{0.1}$ &           Inf &              Inf \\
			st\_e17 &      2 &        $0.18$ &           $19.6$ &          $1.92$ &            $13.9$ & \textbf{GOpt} &  $\mathbf{0.1}$ &       $24.16$ &              $1$ \\
			st\_e19 &      2 & \textbf{GOpt} &           $20.5$ &   \textbf{GOpt} &            $10.7$ & \textbf{GOpt} &  $\mathbf{0.2}$ & \textbf{GOpt} &            $3.4$ \\
			st\_e09 &      2 & \textbf{GOpt} &           $35.9$ &  \textbf{GOpt*} &             $0.2$ & \textbf{GOpt} &  $\mathbf{0.1}$ &        $0.14$ &            $3.7$ \\
			st\_e23 &      2 & \textbf{GOpt} &           $14.1$ &  \textbf{GOpt*} &             $0.3$ & \textbf{GOpt} &  $\mathbf{0.1}$ &         $0.2$ &            $3.2$ \\
			st\_e26 &      2 & \textbf{GOpt} &           $15.5$ &  \textbf{GOpt*} &             $0.2$ & \textbf{GOpt} &  $\mathbf{0.1}$ & \textbf{GOpt} &            $2.8$ \\
			ex4\_1\_8 &      2 & \textbf{GOpt} &           $21.6$ &   \textbf{GOpt} &            $47.8$ & \textbf{GOpt} &  $\mathbf{0.1}$ & \textbf{GOpt} &            $3.1$ \\
			st\_e22 &      2 & \textbf{GOpt} &           $15.7$ &  \textbf{GOpt*} &             $0.4$ & \textbf{GOpt} &  $\mathbf{0.1}$ & \textbf{GOpt} &            $2.8$ \\
			st\_e01 &      2 & \textbf{GOpt} &           $14.8$ &  \textbf{GOpt*} &             $0.1$ & \textbf{GOpt} &    $\mathbf{0}$ & \textbf{GOpt} &            $0.9$ \\
			st\_e08 &      2 & \textbf{GOpt} &           $27.6$ &  \textbf{GOpt*} &             $0.3$ & \textbf{GOpt} &  $\mathbf{0.1}$ & \textbf{GOpt} &            $1.8$ \\
			ex4\_1\_9 &      2 & \textbf{GOpt} &           $26.3$ &   \textbf{GOpt} &            $23.7$ & \textbf{GOpt} &  $\mathbf{0.1}$ & \textbf{GOpt} &            $1.8$ \\
			st\_e24 &      2 & \textbf{GOpt} &           $10.8$ &  \textbf{GOpt*} &             $0.3$ & \textbf{GOpt} &  $\mathbf{0.1}$ & \textbf{GOpt} &            $3.2$ \\
			st\_e18 &      2 & \textbf{GOpt} &           $27.1$ &  \textbf{GOpt*} &             $0.5$ & \textbf{GOpt} &  $\mathbf{0.1}$ &           Inf &              Inf \\
			st\_ht &      2 & \textbf{GOpt} &           $13.1$ &  \textbf{GOpt*} &             $0.2$ & \textbf{GOpt} &  $\mathbf{0.1}$ & \textbf{GOpt} &            $3.1$ \\
			st\_cqpjk2 &      3 & \textbf{GOpt} &           $28.3$ &  \textbf{GOpt*} &             $0.2$ & \textbf{GOpt} &  $\mathbf{0.1}$ & \textbf{GOpt} &            $6.3$ \\
			st\_e11 &      3 & \textbf{GOpt} &          $104.5$ &   \textbf{GOpt} &             $5.6$ & \textbf{GOpt} &  $\mathbf{0.1}$ &           Inf &              Inf \\
			st\_e02 &      3 & \textbf{GOpt} &           $42.4$ &  \textbf{GOpt*} &             $0.2$ & \textbf{GOpt} &  $\mathbf{0.1}$ &           Inf &              Inf \\
			ex6\_2\_8 &      3 &       $25.22$ &           $12.2$ &        $100.02$ &            $14.8$ & \textbf{GOpt} &  $\mathbf{2.6}$ &      $100.45$ &            $6.2$ \\
			ex6\_2\_12 &      4 &       $14.82$ &             $11$ &          $0.63$ &            $11.3$ & \textbf{GOpt} &  $\mathbf{4.7}$ &          $17$ &           $13.1$ \\
			ex6\_2\_9 &      4 &        $0.46$ &           $12.5$ &          $0.46$ &    $\mathbf{9.7}$ & \textbf{GOpt} &        $1501.5$ &       $66.21$ &           $13.2$ \\
			st\_bpv1 &      4 & \textbf{GOpt} &           $12.8$ &  \textbf{GOpt*} &             $0.2$ & \textbf{GOpt} &  $\mathbf{0.1}$ & \textbf{GOpt} &           $13.1$ \\
			sample &      4 & \textbf{GOpt} &           $47.9$ &   \textbf{GOpt} &            $37.3$ & \textbf{GOpt} &  $\mathbf{0.3}$ &           Inf &              Inf \\
			st\_e04 &      4 &        $2.77$ &           $71.2$ &   \textbf{GOpt} &            $47.4$ & \textbf{GOpt} &  $\mathbf{0.2}$ &           Inf &              Inf \\
			st\_bpv2 &      4 & \textbf{GOpt} &           $19.4$ &  \textbf{GOpt*} &             $0.3$ & \textbf{GOpt} &  $\mathbf{0.1}$ & \textbf{GOpt} &            $6.5$ \\
			st\_e41 &      4 &        $0.17$ &           $37.8$ &          $0.17$ &            $34.5$ & \textbf{GOpt} &  $\mathbf{0.1}$ & \textbf{GOpt} &           $12.7$ \\
			st\_e12 &      4 & \textbf{GOpt} &            $8.8$ &   \textbf{GOpt} &             $8.7$ & \textbf{GOpt} &  $\mathbf{0.1}$ & \textbf{GOpt} &            $6.1$ \\
			ex6\_2\_14 &      4 & \textbf{GOpt} &           $22.9$ &   \textbf{GOpt} &   $\mathbf{10.2}$ &           Inf &             Inf &           Inf &              Inf \\
			ex3\_1\_2 &      5 &           Inf &              Inf &  \textbf{GOpt*} &             $0.4$ & \textbf{GOpt} &  $\mathbf{0.1}$ & \textbf{GOpt} &           $12.4$ \\
			st\_e05 &      5 & \textbf{GOpt} &           $36.1$ &  \textbf{GOpt*} &             $0.5$ & \textbf{GOpt} &  $\mathbf{0.1}$ &           Inf &              Inf \\
			ex2\_1\_1 &      5 &       $17.65$ &           $25.5$ &  \textbf{GOpt*} &    $\mathbf{0.3}$ & \textbf{GOpt} &        $1501.4$ &       $67.65$ &           $82.7$ \\
			ex6\_2\_13 &      6 & \textbf{GOpt} &           $28.6$ &   \textbf{GOpt} &   $\mathbf{11.1}$ & \textbf{GOpt} &        $1500.8$ & \textbf{GOpt} &           $34.2$ \\
			ex7\_2\_2 &      6 & \textbf{GOpt} &          $161.8$ &   \textbf{GOpt} &             $8.1$ & \textbf{GOpt} &  $\mathbf{0.3}$ &           Inf &              Inf \\
			st\_e21 &      6 & \textbf{GOpt} &           $13.2$ &   \textbf{GOpt} &            $10.9$ & \textbf{GOpt} &  $\mathbf{0.1}$ & \textbf{GOpt} &           $20.7$ \\
			st\_bsj4 &      6 &        $1.54$ &           $54.1$ &  \textbf{GOpt*} &             $0.2$ & \textbf{GOpt} &  $\mathbf{0.1}$ & \textbf{GOpt} &           $45.9$ \\
			ex6\_2\_10 &      6 & \textbf{GOpt} &           $31.6$ &   \textbf{GOpt} &     $\mathbf{13}$ &           Inf &             Inf &           Inf &              Inf \\
			st\_bsj3 &      6 & \textbf{GOpt} &           $61.7$ &  \textbf{GOpt*} &    $\mathbf{0.2}$ & \textbf{GOpt} &  $\mathbf{0.2}$ & \textbf{GOpt} &             $35$ \\
			ex5\_2\_4 &      7 & \textbf{GOpt} &           $71.5$ &  \textbf{GOpt*} &    $\mathbf{0.4}$ & \textbf{GOpt} &  $\mathbf{0.4}$ & \textbf{GOpt} &           $45.4$ \\
			ex3\_1\_1 &      8 & \textbf{GOpt} &          $118.5$ &  \textbf{GOpt*} &   $\mathbf{47.5}$ & \textbf{GOpt} &         $256.1$ &           Inf &              Inf \\
			ex7\_2\_4 &      8 &           Inf &              Inf &          $0.15$ &            $34.3$ & \textbf{GOpt} & $\mathbf{25.2}$ &           Inf &              Inf \\
			ex7\_2\_3 &      8 & \textbf{GOpt} &          $117.1$ &   \textbf{GOpt} &            $56.7$ & \textbf{GOpt} &        $1501.5$ & \textbf{GOpt} &   $\mathbf{4.1}$ \\
			ex5\_4\_2 &      8 & \textbf{GOpt} &          $106.1$ &  \textbf{GOpt*} &               $1$ & \textbf{GOpt} &  $\mathbf{0.3}$ &           Inf &              Inf \\
			st\_iqpbk2 &      8 & \textbf{GOpt} &           $77.4$ &  \textbf{GOpt*} &             $0.5$ & \textbf{GOpt} &  $\mathbf{0.3}$ & \textbf{GOpt} &           $63.3$ \\
			st\_iqpbk1 &      8 & \textbf{GOpt} &             $81$ &  \textbf{GOpt*} &             $0.7$ & \textbf{GOpt} &  $\mathbf{0.3}$ & \textbf{GOpt} &           $61.9$ \\
			ex6\_2\_5 &      9 &        $0.28$ &           $55.2$ &          $0.28$ &   $\mathbf{48.9}$ & \textbf{GOpt} &        $1502.2$ &       $59.83$ &           $77.8$ \\
			ex6\_2\_7 &      9 &        $0.55$ &           $46.6$ &          $0.31$ &   $\mathbf{14.8}$ & \textbf{GOpt} &        $1502.2$ &           Inf &              Inf \\
			ex5\_2\_2\_case2 &      9 & \textbf{GOpt} &           $82.2$ &  \textbf{GOpt*} &             $0.4$ & \textbf{GOpt} &  $\mathbf{0.2}$ & \textbf{GOpt} &            $6.5$ \\
			ex5\_2\_2\_case3 &      9 & \textbf{GOpt} &           $69.9$ &  \textbf{GOpt*} &             $0.3$ & \textbf{GOpt} &  $\mathbf{0.1}$ & \textbf{GOpt} &            $6.3$ \\
			ex5\_2\_2\_case1 &      9 & \textbf{GOpt} &           $80.1$ &  \textbf{GOpt*} &             $0.3$ & \textbf{GOpt} &  $\mathbf{0.1}$ & \textbf{GOpt} &            $5.3$ \\
			st\_e33 &      9 & \textbf{GOpt} &            $124$ &  \textbf{GOpt*} &             $0.3$ & \textbf{GOpt} &  $\mathbf{0.1}$ &           Inf &              Inf \\
			st\_jcbpaf2 &     10 & \textbf{GOpt} &             $37$ &  \textbf{GOpt*} &             $0.5$ & \textbf{GOpt} &  $\mathbf{0.4}$ &       $51.14$ &           $65.7$ \\
			ex2\_1\_5 &     10 & \textbf{GOpt} &           $23.5$ &  \textbf{GOpt*} &             $0.8$ & \textbf{GOpt} &  $\mathbf{0.6}$ & \textbf{GOpt} &           $63.6$ \\
			st\_bpaf1b &     10 & \textbf{GOpt} &           $35.3$ &  \textbf{GOpt*} &             $0.8$ & \textbf{GOpt} &  $\mathbf{0.1}$ & \textbf{GOpt} &           $69.8$ \\
			st\_e03 &     10 &           Inf &              Inf &   \textbf{GOpt} &            $83.2$ & \textbf{GOpt} &  $\mathbf{0.4}$ &           Inf &              Inf \\
			ex2\_1\_6 &     10 &       $25.64$ &             $43$ &  \textbf{GOpt*} &             $0.3$ & \textbf{GOpt} &  $\mathbf{0.2}$ &       $24.62$ &          $141.1$ \\
			st\_bpaf1a &     10 & \textbf{GOpt} &           $31.4$ &  \textbf{GOpt*} &             $0.4$ & \textbf{GOpt} &  $\mathbf{0.1}$ & \textbf{GOpt} &           $81.9$ \\
			process &     10 & \textbf{GOpt} &          $221.3$ &   \textbf{GOpt} &           $137.1$ & \textbf{GOpt} &  $\mathbf{0.7}$ &           Inf &              Inf \\
			st\_e07 &     10 & \textbf{GOpt} &          $148.6$ &  \textbf{GOpt*} &             $0.3$ & \textbf{GOpt} &  $\mathbf{0.1}$ &           Inf &              Inf \\
			st\_e16 &     12 &           Inf &              Inf &         $11.29$ &             $8.6$ & \textbf{GOpt} &  $\mathbf{0.2}$ &           Inf &              Inf \\
			st\_e30 &     14 &        $8.22$ &          $389.6$ &  \textbf{GOpt*} &             $0.4$ & \textbf{GOpt} &  $\mathbf{0.2}$ &           Inf &              Inf \\
			alkyl &     14 &        $22.6$ &          $292.5$ &          $3.46$ &           $118.9$ & \textbf{GOpt} &  $\mathbf{0.2}$ &           Inf &              Inf \\
			ex8\_4\_5 &     15 & \textbf{GOpt} &          $334.3$ &          $3.24$ &           $797.9$ &        $0.33$ &  $\mathbf{3.6}$ &           Inf &              Inf \\
			ex5\_4\_3 &     16 &           Inf &              Inf &   \textbf{GOpt} &            $11.4$ & \textbf{GOpt} &  $\mathbf{0.2}$ &           Inf &              Inf \\
			ex5\_3\_2 &     22 & \textbf{GOpt} &          $214.9$ &  \textbf{GOpt*} &             $0.6$ & \textbf{GOpt} &  $\mathbf{0.5}$ &           Inf &              Inf \\
			ex2\_1\_8 &     24 &       $18.36$ &          $259.3$ &  \textbf{GOpt*} &             $0.8$ & \textbf{GOpt} &  $\mathbf{0.4}$ &           Inf &              Inf \\
			ex8\_4\_2 &     24 &        $0.63$ &           $1167$ &   \textbf{GOpt} &          $1427.7$ & \textbf{GOpt} &        $1501.1$ &      $296.78$ & $\mathbf{723.7}$ \\
			ex5\_2\_5 &     32 &           Inf &              Inf &  \textbf{GOpt*} &   $\mathbf{46.4}$ & \textbf{GOpt} &        $1501.5$ &           Inf &              Inf \\
			ex8\_2\_4a &     55 & \textbf{GOpt} &         $1448.4$ &   \textbf{GOpt} &          $1394.5$ & \textbf{GOpt} &  $\mathbf{1.1}$ &           Inf &              Inf \\
			ex8\_3\_9 &     78 &        $0.31$ &          $567.6$ &  \textbf{GOpt*} &  $\mathbf{167.2}$ &        $2.88$ &        $1501.1$ &           Inf &              Inf \\
			ex8\_3\_3 &    110 &           Inf &              Inf &         $19.83$ &    $\mathbf{168}$ & \textbf{GOpt} &        $1500.8$ &           Inf &              Inf \\
			ex8\_3\_14 &    110 &           Inf &              Inf &   \textbf{GOpt} &  $\mathbf{356.8}$ & \textbf{GOpt} &        $1501.4$ &           Inf &              Inf \\
			ex8\_3\_2 &    110 &        $5.55$ & $\mathbf{214.7}$ &          $8.25$ &          $1178.7$ & \textbf{GOpt} &        $1501.4$ &           Inf &              Inf \\
			ex8\_3\_4 &    110 &           Inf &              Inf &          $5.31$ & $\mathbf{1243.1}$ & \textbf{GOpt} &        $1501.1$ &           Inf &              Inf \\
			\bottomrule
			\caption{GoML vs BARON vs OCTHaGOn.}\label{tab:results}\\
		\end{longtable}

}

\subsection{Attributions}
In order to better understand which of our enhancements
are more important in improving the optimality gaps, 
we ran the first version of our methodology (GoML Base)
in each one of the $74$ instances of Table \ref{tab:results}, with and
without various enhancements (i.e., Robustness, OCT Sampling and Relaxations).
In particular, to test the effect of robustness we ran the experiments with $\rho=0$ (corresponding
to no robustness) and $\rho=0.01,\:\rho=0.1$  corresponding to variable
levels of robustness.
Then, for the constraint relaxations, we first experimented with no relaxations
and then with a relaxation coefficient of $\lambda=100$, which represents a case where
relaxations are used.
Finally, we also tested the method with and without OCT sampling.\par 
Then, for each of the enhancements we mentioned, we measured the number of instances
in which the optimality gap improved under the presence of the enhancement
(compared to the case where the enhancement was not present).
The results are shown in Table \ref{tbl:attribution}.

{\renewcommand{\arraystretch}{1.2}
\begin{table}[h]
		\begin{tabular}{|l|c|c|c|}
			\hline
			& Robustness & OCT Sampling & Relaxations \\ \hline
			\# Instances where gap improved & 11         & 14           & 13          \\ \hline
			\# Instances where gap stayed the same & 26         & 52           & 36          \\ \hline
		\end{tabular}
	\caption{Effect of the presence of the various enhancements on optimality gaps.}
	\label{tbl:attribution}
\end{table}
}

We notice that each one of the enhancements improves the gaps
in at least $11$ of the $77$ instances.
At the same time, for many instances, enabling an enhancement
does not  hurt the optimality gap.
However, there are instances where enabling an enhancement
leads to increased optimality gaps.
For example, there are many instances where
robustness introduces infeasibilities into the approximating MIO,
thus leading to increased optimality gaps.
For this reason, in order to obtain the best results from
our method, we run it on a problem with and without the
enhancements, and we then keep the best solution, as we mentioned in
the result section.
We again note that the procedure of rerunning the method for different
hyperparameters does not take a lot of time, since testing different
regimes of robustness and relaxations can be done without retraining the ML approximators
approximators (ML training takes up the vast majority of the solution time).\par 
Finally, given that one of the improvements of GoML over OCTHaGOn is that
GoML takes into consideration a larger family of ML models besides OCT-Hs,
we also attempted to see which types of models are used more often by the framework.
In order to do that, we ran the framework on the instances of Table \ref{tab:results}
for different sets of hyperparameters (as described before) and
we found the hyperparameters that produced the best solution  for
that instance. 
Then for those hyperparameters, we logged the types
of ML models that the framework used for approximating
the constraints and objective. 
Those models were the ones with the best accuracy and coefficient of determination ($R^2$)
for approximating the constraints and objective respectively.\par
\begin{figure}[H]
	\centering
	\includegraphics[width=0.6\textwidth]{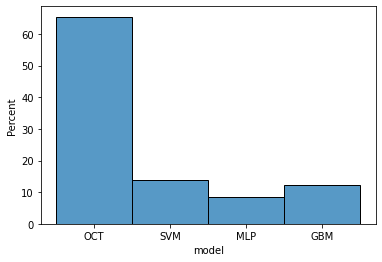}
	\caption{Distribution of models used to approximate constraints/objectives.}
	\label{fig:ml_distr}
\end{figure}

 
The results show that in the MIO approximations, OCT-Hs are used in $65\%$ of the time, SVMs $14\%$ of
the time, GBMs $12\%$ of the time and MLPs $8\%$ of the time, as shown in Figure \ref{fig:ml_distr}.
In particular, OCT-Hs and SVMs tend to be used most of the time in approximating the constraints,
while GBMs and MLPs are used most of the time in approximating the objectives.
Although OCT-Hs are used so often in the constraint approximations, 
they are used very rarely in the objective approximations.
The reason for this is that OCT-Hs require a lot of training
time when trained on a regression task, even for a few variables,
and thus exceed some pre-specified time-limits.
This leads to OCT-Hs not be used that much for approximating
the objectives (i.e., using regression).

\section{Experiments on Real-World and Synthetic Examples}
In this section, we test our approach in some synthetic and real-world
problems.
The purpose of this section is to explore problems
that are not part of the popular MINLP benchmarks of Section \ref{section:experiments}.
\subsection{Quadratic and Sigmoid}
First, we consider the following class of synthetic problems:

\begin{equation}\label{eq:qp_sigm}
	\begin{aligned}
		\min_{\bm{x}} \quad & \bm{c}^T\bm{x}\\
		\textrm{s.t.} \quad &\dfrac{1}{1+{\exp\{-Q_i(\bm{x})\}}}\leq 0.5,\quad i=1,\dots,\lfloor  m/2\rfloor, \\
		&\dfrac{Q_i(\bm{x})}{1+{\exp\{-Q_i(\bm{x})\}}}\geq -0.5,\quad i=\lfloor  m/2\rfloor+1,\dots, m,\\
		&x_k\in[\underbar{x}_k,\overline{x}_k],\: k\in[n],\\
		&\bm{x}\in\RR^n,
	\end{aligned}
\end{equation}
where $Q_i(\bm{x})$ is a quadratic:
\begin{equation}
	Q_i(\bm{x})=\bm{x}^T\bm{A}_i\bm{x}+\bm{d}_i^T\bm{x}+f_i.
\end{equation}
In this problem, the objective is linear
and the constraints involve quadratic and sigmoid functions.
Finally, we also impose bounds on the decision variables,
forcing them to belong to a particular hyper-rectangle.
In order to benchmark our framework in instances of this problem,
we first choose a dimension $n$ and a number of nonlinear constraints $m$.
Then, we randomly generate the cost vector $\bm{c}$ and 
the quadratic parameters $\bm{A}_i,\:\bm{d}_i$ and
$f_i$ for each $i=1,\dots,m$.
We repeat the same process for different values of $n$ and $m$.
Finally, we run the GoML Base framework in the generated instances and
we compare against BARON.\par 
In Table \ref{tbl:sigm} we can see the resulting objective values
returned by the two methods for problem instances of different sizes.
We also record the lower bound calculated by BARON during the solution process.
Finally, since the optimal solution of the instances is not known a-priori, 
we use the lower bound returned by BARON and the best feasible
solution across the two methods to calculate bounds on the optimality
gap of GoML and BARON.
Those bounds are recorded in Table \ref{tbl:sigm}.
\par 
Our experiments show that BARON solves the problem to optimality 
for small values of $n$, but performs poorly for bigger $n$.
For instance, for $n=70$ and $m=2$, BARON terminates early with
a solution that has an optimality gap of at least $99.8\%$,
while the optimality gap of GoML is at most $0.8\%$.
A similar pattern is seen for $n=50$ and $m=4$,
where BARON terminates early with an optimality gap of at least $99.9\%$,
while GoML has an optimality gap of at most $6.56\%$.

{\renewcommand{\arraystretch}{1.2}
\begin{table}[h]
	\begin{tabular}{|l|l|l|l|l|c|c|}
		\hline
		n  & m & GoML Obj. & BARON Obj. & BARON LB & GoML Opt. Gap & BARON Opt. Gap \\ \hline\hline
		$10$ & $2$ & $-22.83$        & $\mathbf{-23.04}$         & -23.04 & $ 0.9\%$  & ${0}$ \\ \hline
		$50$ & $4$ & $\mathbf{-177.27}$        & $-0.225$         & -189.728 & ${\leq  6.56\%}$  & $\geq 99.9\%$ \\ \hline
		$70$ & $2$ & $\mathbf{-255.85}$        & $-0.5711$         & -257.285 & ${\leq 0.8\%}$  & $\geq 99.8\%$ \\ \hline
	\end{tabular}
	\caption{GoML vs BARON in Quadratic-Sigmoid problems.}
	\label{tbl:sigm}
\end{table}
}

\subsection{Speed Reducer Problem}
Following \cite{bertsimas_global_2022}, we also test the method
in the Speed Reducer problem proposed by \cite{jan_golinski_optimal_1970}.
This is a real world problem with the goal of designing a gearbox for an aircraft
engine under geometrical, structural and manufacturing constraints.
The problem is described below:

\begin{equation*}
	\begin{aligned}
		\min _{\mathbf{x}} \quad & 0.7854 x_1 x_2^2\left(3.3333 x_3^2+14.9334 x_3-43.0934\right)  -1.5079 x_1\left(x_6^2+x_7^2\right)+7.477\left(x_6^3+x_7^3\right) \\ 
		\text { s.t. } \quad & -27+x_1 x_2^2 x_3 \geq 0,-397.5+x_1 x_2^2 x_3^2 \geq 0, \\ 
		& -1.93+\frac{x_2 x_6^4 x_3}{x_4^3} \geq 0,-1.93+\frac{x_2 x_7^4 x_3}{x_5^3} \geq 0, \\ 
		& 110.0 x_6^3-\left(\left(\frac{745 x_4}{x_2 x_3}\right)^2+16.9 \times 10^6\right)^{0.5} \geq 0, \\ 
		& 85.0 x_7^3-\left(\left(\frac{745 x_5}{x_2 x_3}\right)^2+157.5 \times 10^6\right)^{0.5} \geq 0, \\ 
		& 40-x_2 x_3 \geq 0, x_1-5 x_2 \geq 0,12 x_2-x_1 \geq 0, \\ & x_4-1.5 x_6-1.9 \geq 0, x_5-1.1 x_7-1.9 \geq 0, \\
		& \mathbf{x} \geq[2.6,0.7,17,7.3,7.3,2.9,5], \\ 
		& \mathbf{x} \leq[3.6,0.8,28,8.3,8.3,3.9,5.5], \\ 
		& x_3 \in \mathbb{Z} .
	\end{aligned}
\end{equation*}
In Table \ref{tbl:srp} we compare the solution given by the enhanced GoML framework
with the solutions given by OCTHaGOn and IPOPT.
We also compare against the optimum by \cite{lin_range_2012}, which according to \cite{bertsimas_global_2022},
is the best known optimum in the literature.
We observe that the ML-based global optimization methods (OCTHaGOn, GoML) perform very well
in this real world problem, being able to find the optimal solution in a reasonable time.
Also, OCTHaGOn, GoML and IPOPT all provide a better solution than that of \cite{lin_range_2012}.\par 
According to the results, IPOPT solves the problem to optimality, beating OCTHaGOn and GoML with respect to solution time.
We note, however, that as described in \cite{bertsimas_global_2022}, in order for IPOPT
to be applied to that problem, the integrality constraint of $x_3$ needs to be relaxed,
since IPOPT does not handle integer variables.
In this particular problem this was not an issue, since the optimal value returned by IPOPT for $x_3$
was indeed an integer.
However, IPOPT cannot be used for general MINLP problems.
On the other hand, the Mixed Integer nature of OCTHaGOn and GoML natively allows
enforcing integrality constraints.
{\renewcommand{\arraystretch}{1.2}
\begin{table}[h]
	\centering
	\begin{tabular}{|r|ccccccc|c|c|}
		\hline & $x_1$ & $x_2$ & $x_3$ & $x_4$ & $x_5$ & $x_6$ & $x_7$ & Objective & Time (s)  \\\hline
		\hline Lin and Tsai (2012) & 3.5 & 0.7 & 17 & 7.3 & 7.7153 & 3.3503 & 5.2867 & 2994.47 & 476 \\
		OCT-HaGOn & 3.5 & 0.7 & 17 & 7.3 & 7.7153 & 3.3502 & 5.2867 & 2994.36 & 32.6 \\
		GoML & 3.5 & 0.7 & 17 & 7.3 & 7.7153 & 3.3502 & 5.2867 & 2994.36 & 37\\
		IPOPT & 3.5 & 0.7 & $17.0^*$ & 7.3 & 7.7153 & 3.3502 & 5.2867 & 2994.36 & 4.2  \\
		\hline
	\end{tabular}
	\caption{Comparison of methods in the speed reducer problem.}	
	\label{tbl:srp}
\end{table}
}
\section{Discussion}
In this section, we discuss the limitations of the approach
and suggest avenues for future work.
\subsection{Contributions}
In this work, we provided extensions to the global optimization
framework OCTHaGOn, in order to improve feasibility and optimality
of the generated solutions.
We tested the framework on a mix of $81$ Global Optimization instances, with $77$ of those
being part the standard benchmarking library MINLPLib and $4$ being
real-world and synthetic instances.
Our results showed that in the majority of instances, our enhancements 
improve the optimality gaps and solutions times of OCTHaGOn.
We also identified $11$ instances where the enhanced framework 
provides better or faster solutions than BARON.\par  
Overall, we showed that despite its approximate nature, 
the enhanced framework, is a promising method
in finding globally optimal solutions in various types of problems.
The framework can potentially be applied in very general global optimization problems,
including problems with convex, non-convex
and even implicit and black-box constraints, with the only assumption
that the user needs to specify bounds for the decision variables.
Hence, due to its generality the method can be used in problems that
are incompatible with traditional global optimizers such as BARON and ANTIGONE.

\subsection{Limitations}
The enhanced framework demonstrates promise in tackling a variety
of global optimization problems, but it is still a work in progress
and comes with its own sets of limitations.\par 
Since the method relies on ML-based approximations
of the original problem, it does not necessarily produce 
globally optimal solutions.
The PGD step at the end of the method helps find
high-quality solutions that are locally optimal,
but without offering guarantees of global optimality
like the ones provided by optimizers such as BARON.
This means that the framework may produce solutions
that are suboptimal or even infeasible in certain cases,
despite the fact that it seems to perform well in a wide range of problems.
At the same time, although the method is designed to handle 
both explicit and implicit constraints, it has only be tested
on explicit problems.
The reason for this is the limited availability of benchmarks that involve
black-box constraints, due to not being solvable by traditional
global optimizers.\par 
Another limitation of the method has to do with the assumptions it makes.
The framework assumes that the decision variables
involved in nonlinear constraints have prespecified bounds,
where in practice bounds may not be available.
In case where bounds are not specified, the method attempts
to compute them through an optimization process, 
although this is not always effective,
since user-specified bounds usually lead to much more precise solutions.
The method also assumes that the nonlinear constraints, 
despite potentially being black-box,
are fast to evaluate.
If this is not the case, such as in simulation-based
black-box constraints,
then the implementation may need to adapt to the computational
needs of the particular application.
Additionally, the PGD step assumes that the
nonlinear constraints support automatic differentiation (AD).
Although this is a soft assumption to make, it may not be true
for general black-box constraints, where AD may not be available.
This, however, can be addressed by approximately evaluating
the gradients (e.g. through finite differencing).\par 
Another consideration about the method is computational speed.
Although the framework is relatively fast for small and medium-sized problems,
the computational time can rise significantly with the number
of variables and nonlinear constraints.
The main bottleneck of the method is the training of the ML models,
and particularly of the hyperplane-based decision trees (OCT-H),
which are used in the sampling and approximation phases of the framework.
To partially account for this problem, we restrict the use of OCT-Hs to smaller
problems, although the ML training time still rises linearly with the number
of nonlinear constraints.
Another speed-related consideration has to do with the complexity of the MIO
approximations.
For small and medium-sized problems, the complexity of the MIO
approximation is usually very small compared to the capabilities
of commercial solvers like Gurobi and CPLEX.
However, an increase in the number of variables accompanied
with an increase in model complexity (e.g.
deeper trees and more layers for MLPs) can significantly
affect MIO solution times.

\section{Conclusion}
In this work, we implemented a range of enhancements to
improve the OCTHaGOn \cite{bertsimas_global_2022} global optimization framework.
We used ML-based sampling, Robust Optimization, and other techniques
to improve the optimality and feasibility gaps of the solutions generated
by OCTHaGOn.
We then demonstrated the effect of our enhancements 
through a range of global optimization benchmarks.
We compared the framework against the commercial optimizer
BARON, and we showed improved solution times and optimality gaps in a  
subset of problems. 
More concretely, we showed that in the majority of test instances,
the enhancements improve the optimality gaps and solution times of OCTHaGOn.
We also showed that in $11$ instances (i.e., $9$ MINLP and $2$ synthetic instances),
the enhanced framework yields better or faster solutions than BARON.\par
These results validate the promising nature of the framework in solving a variety
of problems. 
Overall, the method is a general way of addressing global optimization
problems, that is new in the global optimization literature.
It can handle constraints both explicit and implicit, and even
infer constraints directly from data.
It can deal with arbitrary types of constraints,
without imposing requirements on the mathematical
primitives that are used.
The method only requires bounds of the decision variables involved
in the nonlinear constraints.
Although the method is still new, it can potentially 
be used in a number of different applications, especially
in areas where the optimization problems lack explicit
mathematical formulations.


\bibliography{ML_GOP}

\end{document}